\documentclass[11pt]{article} 
\usepackage{amsfonts,amsmath,latexsym,amssymb,mathrsfs,amsthm}
\usepackage{slashbox}
\usepackage{graphicx}
\usepackage{subcaption}
\usepackage{caption}

\let\OLDthebibliography\thebibliography
\renewcommand\thebibliography[1]{
  \OLDthebibliography{#1}
  \setlength{\parskip}{1pt}
  \setlength{\itemsep}{0pt plus 0.3ex}
}

\def\numberlikeadb{\global\def\theequation{\thesection.\arabic{equation}}}
\numberlikeadb
\newtheorem{theorem}{Theorem}[section]

\newtheorem{conjecture}[theorem]{Conjecture}

\usepackage{lscape}
\usepackage{caption}
\usepackage{multirow}
\usepackage{comment}
\allowdisplaybreaks

\usepackage{color} 

\begin{document}

\title{
The Variance-Gamma Distribution: A Review 
}
\author{Adrian Fischer\footnote{D\'epartement de Math\'ematique, Universit\'e Libre de Bruxelles, Adrian.Fischer@ulb.be}, Robert E. Gaunt\footnote{Department of Mathematics, The University of Manchester, robert.gaunt@manchester.ac.uk}\: and Andrey Sarantsev\footnote{Department of Mathematics \& Statistics, University of Nevada, Reno, asarantsev@unr.edu}}

\date{} 
\maketitle

\thispagestyle{empty}

\vspace{-10mm}

\begin{abstract}The variance-gamma (VG) distributions form a four-parameter family which includes as special and limiting cases the normal, gamma and Laplace distributions. Some of the numerous applications include financial modelling and distributional approximation on Wiener space. In this review, we provide an up-to-date account of the basic distributional theory of the VG distribution. Properties covered include probability and cumulative distribution functions, generating functions, moments and cumulants, mode and median, Stein characterisations, representations in terms of other random variables, and a list of related distributions. We also review methods for parameter estimation and  some applications of the VG distribution, including the aforementioned applications to financial modelling and distributional approximation on Wiener space.
\end{abstract}



\section{Introduction}\label{intro}

The variance-gamma (VG) distribution with parameters $r > 0$, $\theta \in \mathbb{R}$, $\sigma >0$, $\mu \in \mathbb{R}$, denoted by $\mathrm{VG}(r,\theta,\sigma,\mu)$, has probability density function (PDF)
\begin{equation}\label{vgdefn}p(x) = \frac{1}{\sigma\sqrt{\pi} \Gamma(r/2)} \mathrm{e}^{\theta (x-\mu)/\sigma^2} \bigg(\frac{|x-\mu|}{2\sqrt{\theta^2 +  \sigma^2}}\bigg)^{\frac{r-1}{2}} K_{\frac{r-1}{2}}\bigg(\frac{\sqrt{\theta^2 + \sigma^2}}{\sigma^2} |x-\mu| \bigg), \quad x\in\mathbb{R}.
\end{equation}
 In the limit $\sigma\rightarrow 0$ the support becomes the interval $(\mu,\infty)$ if $\theta>0$, and is $(-\infty,\mu)$ if $\theta<0$. Here $K_\nu(x)$ is a modified Bessel function of the second kind; a definition and basic properties that are used throughout the paper are given in Appendix \ref{appa}. 
The parametrisation (\ref{vgdefn}) is taken from Gaunt \cite{gaunt vg} and is similar to one given by Finlay and Seneta \cite{fs06}. Alternative parametrisations are given in Bibby and S{\o}rensen \cite{bibby}
and the book of Kotz, Kozubowski and Podg\'orski \cite{kkp01}, in which they refer to the distribution as the generalized (asymmetric) Laplace distribution. Other names for the VG distribution include the Bessel function distribution (McKay \cite{m32}) and the McKay Type II distribution (Holm and Alouini \cite{ha04}).

The PDF (\ref{vgdefn}) was first written down by Pearson, Jefferey and Elderton  \cite{p29} as the exact PDF of the sample covariance for a random sample drawn from a bivariate normal population. Further distributional properties and applications were given by McKay \cite{m32} and Pearson, Stouffer and David \cite{p32}. The flexibility offered by the four parameters and the modified Bessel function of the second kind in the PDF mean that the VG distribution is often well-suited to statistical modelling;
 for example, Sichel \cite{s73} reported that the VG distribution provided an excellent fit when modelling the size of diamonds mined in South West Africa. The VG distribution (in the symmetric case $\theta=0$) was introduced to the financial literature in a seminal work of Madan and Seneta \cite{madan}, and has since become widely used in financial modelling; see, for example, Madan, Carr and Chang \cite{mcc98}, Madan and Milne \cite{mm91} and Seneta \cite{s04}. The VG distribution is a special case of the generalized hyperbolic (GH) distribution 
that is also popular in financial modelling; see, for example, Eberlein and Keller \cite{eberlein}, Eberlein and Prause \cite{eberlein} and Rydberg \cite{r99}. In addition to its suitability in statistical modelling, the VG distribution has a rich distributional theory with special or limiting cases including the normal, gamma and Laplace distributions, and the product of two zero mean normals and the difference of two independent gammas. In virtue of this, starting with Gaunt \cite{gaunt vg} and Eichelsbacher and Th\"{a}le \cite{eichelsbacher}, the VG distribution has recently found application in the probability literature as a limiting distribution, most notably in analysis on Wiener space.



To date, the most comprehensive account of the VG distribution in the literature is given in Chapter 4 of the excellent book of Kotz et al.\ \cite{kkp01}, in which the VG distribution is viewed as a natural generalisation of the classical (asymmetric) Laplace distribution. As the VG distribution is a special case of the GH distribution, distributional properties can be inferred from results for the GH distribution found in, for example, Bibby and S{\o}rensen \cite{bibby} and Hammerstein \cite{hammersteinphd}. The aforementioned references, however, do not contain even some of the most basic distributional properties of the VG distribution, such as formulas for absolute moments or the mode, and the interested researcher is left to search a rather large and difficult to navigate literature
to track down many such results. 

In this review, we fill in a gap in the literature by providing an up-to-date account of the basic distributional theory of the VG distribution. The end result is that many of the most important distributional properties of the VG distribution are now contained in a single reference. Most results are already stated in the literature or can be inferred from the fact that VG distribution is a special case of the GH distribution and appealing to known results for this distribution. For the small number of results that we could not locate in the literature, we provide straightforward and concise derivations. 

This review covers formulas for the PDF (Section \ref{sec2.1}), the cumulative distribution function (Section \ref{sec2.2}), generating functions, infinite divisibility and self-decomposability (Section \ref{sec2.3}), representations in terms of other random variables (Section \ref{sec2.4}), lists of related distributions (Section \ref{sec2.5}), Stein characterisation (Section \ref{sec2.6}), moments and cumulants (Section \ref{sec2.7}), and mode and median (Section \ref{sec2.8}). The review covers some of the most basic and important properties of a probability distribution, but is not comprehensive; for example, we only briefly touch upon multivariate extensions.

In Section \ref{sec3}, we review methods for parameter estimation for the VG distribution. This section provides a concise overview of a literature that is rather large on account of the popularity of the VG distribution in statistical modelling, particularly in mathematical finance. 
In Section \ref{sec4}, we provide applications of the VG distribution in several research domains.
In Section \ref{sec4.1}, we provide examples in which the VG distribution appears as an exact distribution in connection to sample covariances and Wishart matrices. We review the well-known connection to the variance-gamma process (also referred to as Laplace motion; see for example, Kotz et al.\ \cite{kkp01}) and financial modelling in Sections \ref{sec4.2} and \ref{sec4.3}, respectively. We provide an overview of the recent literature concerning
application of the VG distribution to time series modelling in Section \ref{sec4.4}. In Section \ref{sec4.5}, we review some of the recent activity in the probability literature in which the VG distribution has arisen as a limiting distribution.  

\vspace{2mm}

\noindent{\emph{Notation}.} To simplify formulae, we define 
$
\lambda_{\pm}:= (\sqrt{\theta^2+\sigma^2}\pm\theta)/\sigma^2.
$ 



\section{Distributional Properties}\label{sec2}

\subsection{Density and parametrisations}\label{sec2.1}
We begin by recalling some other parametrisations given in the literature. The first is given in Finlay and Seneta \cite{fs06}. For $x \in \mathbb R$,
\begin{align}\label{vgpam2}p(x)=\sqrt{\frac{2}{\pi}}\frac{\alpha^\alpha \mathrm{e}^{\theta_0(x-\mu)/\sigma_0^2}}{\sigma_0\Gamma(\alpha)}\bigg(\frac{|x-\mu|}{\sqrt{\theta_0^2+2\alpha\sigma_0^2}}\bigg)^{\alpha-1/2}K_{\alpha-1/2}\bigg(\frac{\sqrt{\theta_0^2+2\alpha\sigma_0^2}}{\sigma_0^2}|x-\mu|\bigg).
\end{align}
 It is related to the parametrisation in (\ref{vgdefn}) by $r=2\alpha$, $\sigma^2=\sigma_0^2/(2\alpha)$, $\theta=\theta_0/(2\alpha)$. If a random variable has PDF (\ref{vgpam2}), we write $X\sim \mathrm{VG}_2(\alpha,\theta_0,\sigma_0,\mu)$. Setting $\alpha=1/\nu$ gives the parametrisation of Madan et al.\ \cite{mcc98}, and further setting $\theta_0=0$ yields the parametrisation of the (symmeric) VG distribution of Madan and Seneta \cite{madan}. Another parametrisation is given in Bibby and S{\o}rensen \cite{bibby}:
\begin{equation}\label{vgpam3} p(x) = \frac{\gamma^{2\lambda}}{\sqrt{\pi} \Gamma(\lambda)}\left(\frac{|x-\mu|}{2\alpha}\right)^{\lambda-1/2} \mathrm{e}^{\beta (x-\mu)}K_{\lambda-1/2}(\alpha|x-\mu|), \quad x\in \mathbb{R},
\end{equation}
where $\gamma^2=\alpha^2-\beta^2$, and is related to the parametrisation in (\ref{vgdefn}) by $r=2\lambda$, $\theta=\beta/\gamma^2$, $\sigma=1/\gamma$. Finally, we note the following parametrisation used by Kotz et al.\ \cite{kkp01}:
\begin{equation*}p(x)=\frac{\sqrt{2}\mathrm{e}^{\frac{\sqrt{2}}{2\sigma_0}(1/\kappa-\kappa)(x-\mu)}}{\sigma_0^{\tau+1/2}\Gamma(\tau)}\bigg(\frac{\sqrt{2}|x-\mu|}{\kappa+1/\kappa}\bigg)^{\tau-1/2}K_{\tau-1/2}\bigg(\frac{\sqrt{2}}{2\sigma_0}\bigg(\frac{1}{\kappa}-\kappa\bigg)|x-\mu|\bigg), \quad x\in\mathbb{R},
\end{equation*} 
which is related to the parametrisation in (\ref{vgdefn}) by $r=2\tau$, $\theta=\sigma_0(1/\kappa-\kappa)/2^{3/2}$, $\sigma^2=\sigma_0^2/2$. Henceforth, we shall work mostly with the parametrisation (\ref{vgdefn}), with results in the other parametrisations being readily deduced using the above change of parameters.
 
\begin{figure}[!]
    \centering
\vspace{.2cm}
   \begin{subfigure}{.49\textwidth}
\captionsetup{width=.95\textwidth}
  \centering
  \includegraphics[width=7.4cm]{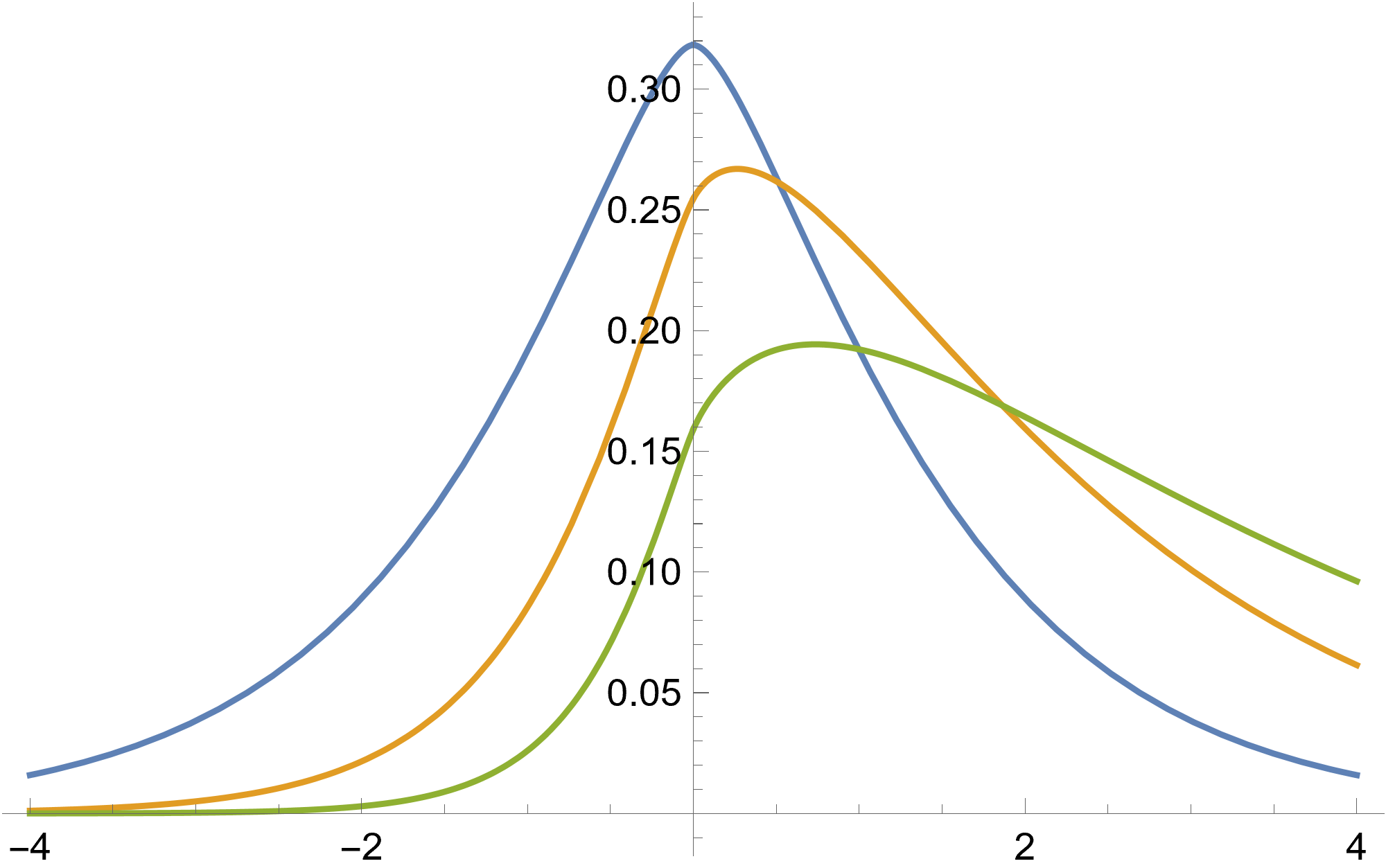}
\end{subfigure}%
  \begin{subfigure}{.49\textwidth}
\captionsetup{width=.95\textwidth}
  \centering
  \includegraphics[width=7.4cm]{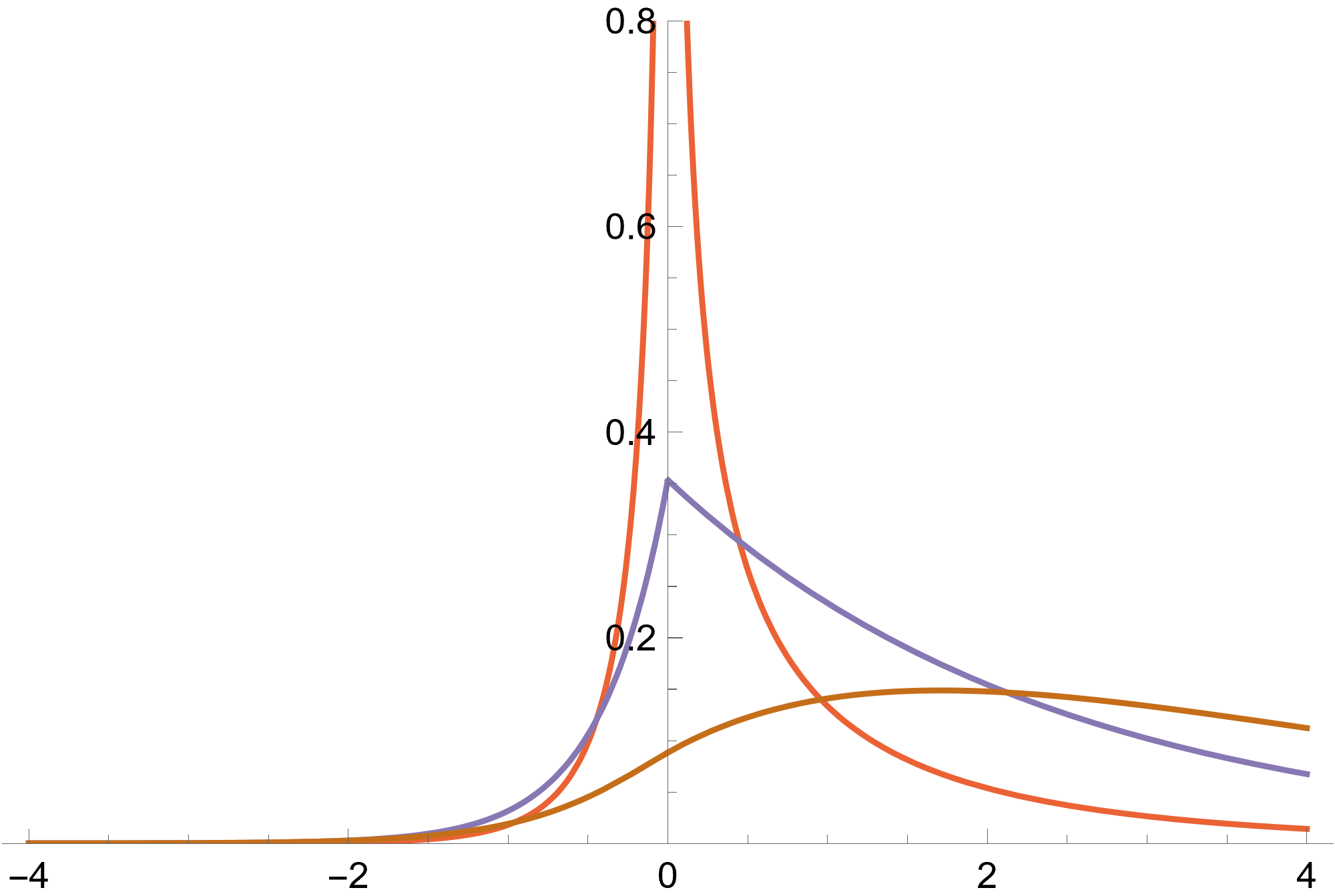}
\end{subfigure}
\caption{\label{fig:pdf}
\it 
The VG PDF for different parameter constellations: \\ \includegraphics[scale=1]{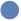}$\mathrm{VG}(3,0,1,0)$, \includegraphics[scale=1]{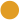}$\mathrm{VG}(3,0.5,1,0)$, \includegraphics[scale=1]{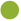}$\mathrm{VG}(3,1,1,0)$ (left image) and \\ \includegraphics[scale=1]{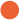}$\mathrm{VG}(0.5,1,1,0)$, \includegraphics[scale=1]{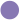}$\mathrm{VG}(2,1,1,0)$, \includegraphics[scale=1]{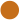}$\mathrm{VG}(4,1,1,0)$ (right image).
} 
\end{figure}

The presence of the modified Bessel function in the PDF (\ref{vgdefn}) makes it a little difficult to parse on first inspection.  The following limiting forms can help in this regard. Suppose $r > 0$, $\theta \in \mathbb{R}$, $\sigma >0$, $\mu \in \mathbb{R}$.  Using the limiting form (\ref{Ktend0}) gives that, as $x\rightarrow\mu$, 
\begin{equation}\label{pmutend}p(x)\sim\begin{cases}\displaystyle \frac{(1+\theta^2/\sigma^2)^{-(r-1)/2}}{2\sigma\sqrt{\pi}}\frac{\Gamma((r-1)/2)}{\Gamma(r/2)}, &  r>1, \\
\displaystyle -\frac{1}{\pi\sigma}\log|x-\mu|, &  r=1, \\
\displaystyle \frac{1}{(2\sigma)^r\sqrt{\pi}}\frac{\Gamma((1-r)/2)}{\Gamma(r/2)}|x-\mu|^{r-1}, &  0<r<1 \end{cases}
\end{equation}
(see Gaunt \cite{gaunt thesis}).  We see from (\ref{pmutend}) that the PDF has a singularity at the origin for $r\leq1$.  Moreover, for all parameter values ($r>0$, $\theta\in\mathbb{R}$, $\sigma>0$, $\mu\in\mathbb{R}$), the $\mathrm{VG}(r,\theta,\sigma,\mu)$ distribution is unimodal; further details are given in Section \ref{sec2.8}. The density is bounded for $r>1$.  
Plots of the VG PDF (\ref{vgdefn}) that show the effect of varying the shape parameter $r$ (which agree with these assertions) and the skewness parameter $\theta$ are given in Figure \ref{fig:pdf}.
Also, applying the limiting form (\ref{Ktendinfinity}) to (\ref{vgdefn}) gives that (see Gaunt \cite{gaunt vg3})
\begin{equation}\label{tail1}
p(x)\sim \frac{x^{r/2-1}}{2^{r/2}(\theta^2+\sigma^2)^{r/4}\Gamma(r/2)}\mathrm{e}^{-\lambda_-(x-\mu)}, \quad x\rightarrow\infty,
\end{equation}
\begin{equation}\label{tail2}
p(x)\sim\frac{(-x)^{r/2-1}}{2^{r/2}(\theta^2+\sigma^2)^{r/4}\Gamma(r/2)}\mathrm{e}^{\lambda_+(x-\mu)}, \quad x\rightarrow-\infty.
\end{equation}
We observe that the tails of the VG distribution are heavier than the tails of the normal distribution.
This feature is important in financial modelling and allows for a better fit to real financial data than the normal distribution; see Section \ref{sec4.3} for further details. 

When $r\in\mathbb{Z}^+=\{1,2,3,\ldots\}$ is even, we can use a standard simplification of the modified Bessel function of the second kind (see (\ref{special})) to obtain a formula for the PDF in terms of elementary functions (see Kotz et al.\ \cite{kkp01}, and earlier Teichroew \cite{t57} for the case $\theta=0$):
\begin{align}p(x)&= \frac{|x-\mu|^{r/2-1}}{2^{r/2}(\theta^2+\sigma^2)^{r/4}\Gamma(r/2)}\exp\bigg(\frac{\theta}{\sigma^2}(x-\mu)-\frac{\sqrt{\theta^2+\sigma^2}}{\sigma^2}|x-\mu|\bigg)\times\nonumber\\
\label{nsimp}&\quad\times\sum_{j=0}^{r/2-1}\frac{(r/2-1+j)!}{(r/2-1-j)!j!}\bigg(\frac{\sigma^2}{2\sqrt{\theta^2+\sigma^2}|x-\mu|}\bigg)^j, \quad x\in\mathbb{R},\:r\in 2\mathbb{Z}^+.
\end{align} 
 
It is immediate from (\ref{vgdefn}) that the class of VG distributions is closed under affine transformations (see Kotz et al.\ \cite{kkp01}). Let $X\sim \mathrm{VG}(r,\theta,\sigma,\mu)$. Then, for $a\not=0$ and $b\in\mathbb{R}$,
\begin{align}\label{muc}aX+b\sim\mathrm{VG}(r,a\theta,|a|\sigma,a\mu+b).
\end{align}
To simplify expressions, we shall sometimes set $\mu=0$, with results for the general case $\mu\in\mathbb{R}$ immediately following because $\mu+\mathrm{VG}(r,\theta,\sigma,0)=_d\mathrm{VG}(r,\theta,\sigma,\mu)$ (with obvious abuse of notation). We also observe that if $X\sim \mathrm{VG}(r,\theta,\sigma,\mu)$ then $-X\sim\mathrm{VG}(r,-\theta,\sigma,-\mu)$, whilst if $X\sim\mathrm{VG}(r,0,\sigma,0)$ then $-X=_d X$. In line with terminology of Kotz et al.\ \cite{kkp01}, we say that $X\sim\mathrm{VG}(r,0,\sigma,0)$ has a \emph{symmetric variance-gamma} distribution.

\subsection{Cumulative distribution function}\label{sec2.2}

A closed-form formula for the cumulative distribution function (CDF) of the VG distribution is not available for general parameter values $r>0$, $\theta\in\mathbb{R}$, $\sigma>0$, $\mu\in\mathbb{R}$. We note some cases for which exact formulas can be given. For $X\sim\mathrm{VG}(r,\theta,\sigma,\mu)$, let $F(x)=\mathbb{P}(X\leq x)$. 

Suppose $\theta=0$. Then, by the symmetry of the $\mathrm{VG}(r,0,\sigma,\mu)$ distribution about $x=\mu$, it follows that $F(x)=1/2+\mathrm{sgn}(x)\int_0^{|x|}p(t)\,\mathrm{d}t$, where $\mathrm{sgn}(x)$ is the sign of $x$. Calculating the integral using (\ref{besint}), gives that, for $x\in\mathbb{R}$,
\begin{align*}F(x)=\frac{1}{2}+\frac{(x-\mu)}{2\sigma}\bigg[K_{\frac{r-1}{2}}\bigg(\frac{|x-\mu|}{\sigma}\bigg)\mathbf{L}_{\frac{r-3}{2}}\bigg(\frac{|x-\mu|}{\sigma}\bigg)+\mathbf{L}_{\frac{r-1}{2}}\bigg(\frac{|x-\mu|}{\sigma}\bigg)K_{\frac{r-3}{2}}\bigg(\frac{|x-\mu|}{\sigma}\bigg)\bigg],
\end{align*}
where $\mathbf{L}_\nu(x)$ is a modified Struve function of the first kind (see Olver et al.\ \cite[Chapter 11]{olver} for a definition and basic properties).
Other formulas for the CDF for the case $\theta=0$ are given by Jankov Ma\v{s}irevi\'c and Pog\'any \cite{jp} and Nadarajah, Srivastava and Gupta \cite{nsg07}.

Now suppose $r\in 2\mathbb{Z}^+$ and $\theta\in\mathbb{R}$. Then, making use of the formula (\ref{nsimp})
one readily obtains the following formulas (see Nadarajah et al.\ \cite{nsg07}). For $x\leq\mu$,
\begin{align*}F(x)&=\frac{(\theta^2+\sigma^2)^{-r/4}}{(2\lambda_+)^{r/2}(r/2-1)!}\sum_{j=0}^{r/2-1}\frac{(r/2-1+j)!}{(r/2-1-j)!j!}\bigg(\frac{\theta+\sqrt{\theta^2+\sigma^2}}{\sqrt{\theta^2+\sigma^2}}\bigg)^j\Gamma\bigg(\frac{r}{2}-j,-\lambda_+(x-\mu)\bigg),
\end{align*}
and, for $x>\mu$,
\begin{align*}F(x)&=1-\frac{(\theta^2+\sigma^2)^{-r/4}}{(2\lambda_-)^{r/2}(r/2-1)!}\sum_{j=0}^{r/2-1}\frac{(r/2-1+j)!}{(r/2-1-j)!j!}\bigg(\frac{\sqrt{\theta^2+\sigma^2}-\theta}{\sqrt{\theta^2+\sigma^2}}\bigg)^j\Gamma\bigg(\frac{r}{2}-j,\lambda_-(x-\mu)\bigg),
\end{align*}
where $\Gamma(a,x)=\int_x^\infty t^{a-1}\mathrm{e}^{-t}\,\mathrm{d}t$ is the upper incomplete gamma function. 

Let us now return to the general setting $r>0$, $\theta\in\mathbb{R}$, $\sigma>0$, $\mu\in\mathbb{R}$. Since a closed-form formula is not available for the CDF, asymptotic approximations for the tail probabilities are of interest. Let $\bar{F}(x)=1-F(x)=\mathbb{P}(X> x)$ for $X\sim\mathrm{VG}(r,\theta,\sigma,\mu)$. Then, using (\ref{vgdefn}) and the limiting form (\ref{kintap}) gives that, as $x\rightarrow\infty$,
\begin{align}\label{barbar}\bar{F}(x)\sim\frac{1}{2^{r/2}(\theta^2+\sigma^2)^{r/4}\lambda_-\Gamma(r/2)}x^{r/2-1}\mathrm{e}^{-\lambda_-(x-\mu)}, 
\end{align}
and, by symmetry, as $x\rightarrow-\infty$,
\begin{align*}F(x)\sim\frac{1}{2^{r/2}(\theta^2+\sigma^2)^{r/4}\lambda_+\Gamma(r/2)}(-x)^{r/2-1}\mathrm{e}^{\lambda_+(x-\mu)}.
\end{align*}
 
Upper and lower bounds for $F(x)$ and $\overline{F}(x)$ can be obtained by using bounds of Gaunt \cite{gaunt ineq1, gaunt ineq3} for the integral $\int_x^\infty\mathrm{e}^{\beta t}t^\nu K_{\nu}(t)\,\mathrm{d}t$, $x>0$, $-1<\beta<1$, $\nu>-1/2$. As an example, from inequality (2.9) of Gaunt \cite{gaunt ineq3} we have that $\int_x^\infty\mathrm{e}^{\beta t}t^\nu K_{\nu}(\alpha t)\,\mathrm{d}t\leq (\alpha-\beta)^{-1}\mathrm{e}^{\beta x}x^\nu K_\nu(\alpha x)$, for $x>0$, $0\leq\beta<\alpha$, $-1/2<\nu\leq1/2$, with equality if and only if $\nu=1/2$, and the inequality is reversed if $\nu>1/2$. Applying this bound to (\ref{vgdefn}) yields that, for $\theta\geq0$, $x>\mu$,
\begin{align}
\label{cumbd}\bar{F}(x)&\leq \frac{ \mathrm{e}^{\theta (x-\mu)/\sigma^2} }{\sigma\sqrt{\pi} \lambda_-\Gamma(r/2)} \bigg(\frac{|x-\mu|}{2\sqrt{\theta^2 +  \sigma^2}}\bigg)^{\frac{r-1}{2}} K_{\frac{r-1}{2}}\bigg(\frac{\sqrt{\theta^2 + \sigma^2}}{\sigma^2} |x-\mu| \bigg)=\frac{p(x)}{\lambda_-},
\end{align} 
with equality if and only if $r=2$, and the inequality is reversed if $r>2$. Here $p$ is the PDF (\ref{vgdefn}). Inequality (\ref{cumbd}) is tight as $x\rightarrow\infty$, which can be seen by applying the limiting form (\ref{Ktendinfinity}) to the bound in (\ref{cumbd}) and comparing to (\ref{barbar}).

\subsection{Generating functions, infinite divisibility and self-decomposability}\label{sec2.3}

The moment generating function of $X\sim\mathrm{VG}(r,\theta,\sigma,\mu)$ is easily calculated using the integral formula (\ref{pdfintegral}) (see, for example, Bibby and S{\o}rensen \cite{bibby}),
\begin{align}\label{mgf}M(t)=\mathbb{E}[\mathrm{e}^{tX}]=\mathrm{e}^{\mu t}\big(1-2\theta t-\sigma^2t^2\big)^{-r/2},
\end{align}
which exists provided $-\lambda_+<t<\lambda_-$ (see (\ref{tail1}) and (\ref{tail2})). The characteristic function is
\begin{equation}\label{cfcf}\varphi(t)=\mathbb{E}[\mathrm{e}^{\mathrm{i}tX}]=\mathrm{e}^{\mathrm{i}\mu t}\big(1-2\mathrm{i}\theta t+\sigma^2t^2\big)^{-r/2}.
\end{equation}
The cumulant generating function is defined for $-\lambda_+<t<\lambda_-$, and given by
\begin{align}\label{kfkf}K(t)=\log\mathbb{E}[\mathrm{e}^{tX}]&=\mu t-\frac{r}{2}\log\big(1-2\theta t-\sigma^2t^2\big) \nonumber \\
&=\mu t-\frac{r}{2}\log\bigg(1-\frac{t}{\lambda_-}\bigg)-\frac{r}{2}\log\bigg(1+\frac{t}{\lambda_+}\bigg).
\end{align}

From the characteristic function (\ref{cfcf}) formula, it is clear that if $X_1\sim\mathrm{VG}(r_1,\theta,\sigma,\mu_1)$ and $X_2\sim\mathrm{VG}(r_2,\theta,\sigma,\mu_2)$ are independent, then
\begin{equation}\label{conv}X_1 +X_2 \sim \mathrm{VG}(r_1+r_2 ,\theta,\sigma,\mu_1+\mu_2)
\end{equation}
(see Bibby and S{\o}rensen \cite{bibby}). Thus, the class of VG distributions is closed under convolution, provided the random variables have common values of $\theta$ and $\sigma$. 

It is clear from (\ref{conv}) that the $\mathrm{VG}(r,\theta,\sigma,\mu)$ distribution is infinitely divisible.
Moreover, the VG characteristic function has the following L\'{e}vy-Hinchin representation (see Madan et al.\ \cite{mcc98}):
\begin{align}\label{levyrep}\varphi(t)=\exp\bigg(\mathrm{i}\mu t+\int_{-\infty}^\infty(\mathrm{e}^{\mathrm{i}tx}-1
)\nu(x)\,\mathrm{d}x \bigg),
\end{align}
where the L\'evy density is given by
\begin{equation}\label{levymeas}\nu(x)=\frac{r}{2|x|}\mathrm{e}^{-\lambda_+|x|}\mathbf{1}_{x<0}+\frac{r}{2x}\mathrm{e}^{-\lambda_-x}\mathbf{1}_{x>0}.
\end{equation}

A distribution $\mathcal Q$ on the real line is {\it self-decomposable}  (Sato \cite[Definition 15.1, p.\ 90]{sato99}) if for any constant $c \in (0, 1)$ there exists a probability space with a random variable $X \sim \mathcal Q$ and independent
random variable $Y$ such that $cX + Y \sim \mathcal Q$. 
Self-decomposable distributions are infinitely divisible and have a special form of L\'evy-Hinchin representation. As the L\'evy density $\nu(x)$, given by (\ref{levymeas}), can be written in the form $\nu(x)=h(x)/|x|$, where $h(x)\geq0$ is increasing for negative $x$ and decreasing for positive $x$, it is immediate from the representation (\ref{levyrep}) (see Sato \cite[Corollary 15.11]{sato99}) that the $\mathrm{VG}(r,\theta,\sigma,\mu)$ distribution is self-decomposable. An explicit construction of the self-decomposability of the $\mathrm{VG}(r,\theta,\sigma,\mu)$ distribution is also available for the case that $r\in2\mathbb{Z}^+$. Let $X\sim\mathrm{VG}(2n,\theta,\sigma,\mu)$, where $n\in\mathbb{Z}^+$. Then, for any $c\in[0,1]$,
\begin{align}\label{vgselfdec}X=_d cX+(1-c)\mu+\sum_{i=1}^nV_i,
\end{align}
where $V_1,\ldots,V_n$ are i.i.d.\ random variables that can be expressed as $V_1=_d \sigma(\delta_1W_1/(\sigma\lambda_-)-\sigma\lambda_-\delta_2W_2)$. Here, $\delta_1$ and $\delta_2$ are $0$-$1$ random variables 
with probabilities
\begin{align*}&\mathbb{P}(\delta_1=0,\delta_2=0)=c^2, \quad 
\mathbb{P}(\delta_1=1,\delta_2=0)=(1-c)\bigg(c+\frac{1-c}{1+\sigma^2\lambda_-^2}\bigg), \\
& \mathbb{P}(\delta_1=0,\delta_2=1)=(1-c)\bigg(c+\frac{(1-c)\sigma^2\lambda_-^2}{1+\sigma^2\lambda_-^2}\bigg), \quad \mathbb{P}(\delta_1=1,\delta_2=1)=0,
\end{align*}
where $W_1$ and $W_2$ are exponential with rate parameter 1, and $X$, $W_1$, $W_2$, and $(\delta_1, \delta_2)$ are mutually independent. The representation (\ref{vgselfdec}) is given in Kotz et al.\ \cite[Proposition 4.1.4]{kkp01}, in which they generalised the representation of Ramachandran \cite{r97} for the asymmetric Laplace distribution by combining the convolution property (\ref{conv}) of the VG distribution together with the fact that the $\mathrm{VG}(2,\theta,\sigma,\mu)$ distribution is an asymmetric Laplace distribution (see part 1 of Section \ref{subcalss}).

\subsection{Representation in terms of other random variables}\label{sec2.4}

\noindent{1.} The VG distribution has the following fundamental representation in terms of independent normal and gamma random variables; see Kotz et al.\ \cite[Proposition 4.1.2]{kkp01} for a statement of the result and a short proof involving characteristic functions. To fix notation, consider the $\Gamma(r,\lambda)$ distribution with PDF $p(x)=\lambda^rx^{r-1}\mathrm{e}^{-\lambda x}/\Gamma(r)$, $x>0$.
Suppose that  $S\sim \Gamma(r/2,1/2)$ and $T\sim N(0,1)$ are independent random variables.
Then
\begin{equation}\label{rep2}
\mu + \theta S +\sigma \sqrt{S} T\sim \mathrm{VG}(r,\theta,\sigma,\mu).
\end{equation}
The VG distribution is therefore a univariate normal variance-mean distribution (Barndorff-Nielsen, Kent and S{\o}rensen \cite{bks82}). Indeed, conditional on $S$, the random variable $\mu + \theta S +\sigma \sqrt{S} T$ has the $N(\mu+\theta S, \sigma^2S)$ distribution. 

\vspace{2mm}

\noindent{2.} When $r\geq1$ is a positive integer, the VG distribution can be represented in terms of independent standard normal random variables $X_1,\ldots,X_r$ and $Y_1,\ldots,Y_r$:
\begin{equation}\label{gauntrep}\mu + \theta \sum_{i=1}^r X_i^2 + \sigma \sum_{i=1}^r X_i Y_i\sim \mathrm{VG}(r,\theta,\sigma,\mu). 
\end{equation}
(see Gaunt \cite[Corollary 2.5]{gaunt vg}). To see this, first we define $Z_1=\mu+\theta X_1^2 +\sigma X_1Y_1$ and $Z_i = \theta X_i^2 +\sigma X_iY_i$, $i=2,\ldots,r$. Observe that, for $i=1,2,\ldots,r$, $X_iY_i=_d|X_i|Y_i$ and that $X_i^2\sim\Gamma(r/2,1/2)$.  Hence, by (\ref{rep2}), we have that $Z_1\sim\mathrm{VG}_1(1,\theta,\sigma,\mu)$ and $Z_i\sim\mathrm{VG}_1(1,\theta,\sigma,0)$, $i=2,\ldots, r$.  It therefore follows from (\ref{conv}) that $\sum_{i=1}^rZ_i\sim\mathrm{VG}(r,\theta,\sigma,\mu)$. 

\vspace{2mm}

\noindent{3.} The VG distribution has a neat representation as a difference of independent gamma random variables (this was shown in the $\theta=0$ case by McLeish \cite{m82}; for the general case see, for example, Press \cite{p67} or Kotz et al.\ \cite{kkp01}). Suppose $S\sim\Gamma(r/2,(\sqrt{\theta^2+\sigma^2}+\theta)^{-1})$ and $S'\sim\Gamma(r/2,(\sqrt{\theta^2+\sigma^2}-\theta)^{-1})$ are independent. Then
\begin{equation}\label{gamrep}\mu+S-S'\sim\mathrm{VG}(r,\theta,\sigma,\mu).
\end{equation}
This representation is efficiently proved using a standard characteristic function argument with the $\mathrm{VG}(r,\theta,\sigma,\mu)$ characteristic function formula (\ref{cfcf}) and the formula for the characteristic function of the $\Gamma(r,\lambda)$ distribution, $\varphi_S(t)=(1-\mathrm{i}t/\lambda)^{-r}$.

Suppose now that $r\geq1$ is an integer. Then $S=_d (\sqrt{\theta^2+\sigma^2}+\theta) V/2$ and $S'=_d (\sqrt{\theta^2+\sigma^2}-\theta) V'/2$, where $V$ and $V'$ are independent $\chi_r^2$ random variables (chi-square distribution with $r$ degrees of freedom). From the representation of the $\chi_r^2$ distribution as the sum of the squares of $r$ independent $N(0,1)$ random variables and the representation (\ref{gamrep}), we deduce that 
\[\mu+r\theta+\sum_{j=1}^{2r}\lambda_j(N_j^2-1)\sim \mathrm{VG}(r,\theta,\sigma,\mu),\]
where $N_1,\ldots,N_{2r}$ are independent $N(0,1)$ random variables,  $\lambda_1=\ldots=\lambda_r=(\sqrt{\theta^2+\sigma^2}+\theta)/2$ and $\lambda_{r+1}=\ldots=\lambda_{2r}=-(\sqrt{\theta^2+\sigma^2}-\theta)/2$. So, for integer $r\geq1$ the VG distribution is a member of the second Wiener chaos (see Nourdin and Peccati \cite[Section 2.2]{np12}).

\vspace{2mm}

\noindent{4.} If $r\in 2\mathbb{Z}^+$, the independent random variables $V\sim\chi_r^2$ and $V'\sim\chi_r^2$ can be represented in terms of independent uniform $U(0,1)$ random variables $U_1,\ldots,U_{r}$ by $V=_d\sum_{j=1}^{r/2}\log U_j$ and $V'=_d\sum_{j=r/2+1}^{r}\log U_j$. Thus, from (\ref{gamrep}), 
\begin{equation*}\mu+\frac{\sqrt{\theta^2+\sigma^2}+\theta}{2}\sum_{j=1}^{r/2}\log U_j-\frac{\sqrt{\theta^2+\sigma^2}-\theta}{2}\sum_{j=r/2+1}^{r}\log U_j\sim \mathrm{VG}(r,\theta,\sigma,\mu).
\end{equation*}
This representation is convenient for simulating the VG distribution when $r\in 2\mathbb{Z}^+$. For other values of $r>0$, the $\mathrm{VG}(r,\theta,\sigma,\mu)$ distribution can be simulated using the representation (\ref{gamrep}) and simulating the gamma distributions of $V$ and $V'$ using the methods given in Chapter 9, Section 3 of Devroye \cite{d86}. 

\subsection{Related distributions}\label{sec2.5}

\subsubsection{Subclasses and limiting cases}\label{subcalss}

\noindent{1.} The $\mathrm{VG}(2,\theta,\sigma,\mu)$ distribution corresponds to the asymmetric Laplace distribution. Indeed, setting $r=2$ in (\ref{nsimp}) yields
\begin{equation}\label{lapd}p(x)=\frac{1}{\sqrt{\theta^2+\sigma^2}}\exp\bigg(\frac{\theta}{\sigma^2}(x-\mu)-\frac{\sqrt{\theta^2+\sigma^2}}{\sigma^2}|x-\mu|\bigg).
\end{equation}
Setting $\theta=0$ in (\ref{lapd}) gives the PDF of the classical Laplace distribution. A comprehensive account of the distributional theory of the Laplace and asymmetric Laplace distributions is given in
Chapters 2 and 3 of Kotz et al.\ \cite{kkp01}.

\vspace{2mm}

\noindent{2.} The gamma distribution is a limiting case of the VG distribution (see, for example, Kotz et al.\ \cite{kkp01}). For fixed, $r,\lambda>0$, the sequence of random variables $X_\sigma\sim\mathrm{VG}(2r,(2\lambda)^{-1},\sigma,0)$ converges in distribution to a $\Gamma(r,\lambda)$ random variable, as $\sigma\rightarrow0$. This can be seen by letting $\sigma\rightarrow0$ in the formula (\ref{cfcf}) for the characteristic function of VG distribution and comparing to the gamma characteristic function using L\'evy's continuity theorem. 

\vspace{2mm}

\noindent{3.} The normal distribution is also a limiting case of the VG distribution (see, for example, Kotz et al.\ \cite{kkp01}).  For fixed, $\theta\in\mathbb{R}$ and $\sigma>0$, the sequence of random variables $Y_r\sim\mathrm{VG}(r,\theta/\sqrt{r},\sigma/\sqrt{r},-\theta\sqrt{r})$ converges to a $N(0,\sigma^2+2\theta^2)$ random variable, as $r\rightarrow\infty$.  As in part 2, this can be deduced from a characteristic function argument. More insightfully, from the convolution property, for  $r\in\mathbb{Z}^+$, we can write $Y_r=_d r^{-1/2}\sum_{i=1}^rY_{r,i}$, where $Y_{r,i}\sim\mathrm{VG}(1,\theta,\sigma,-\theta)$. The result now follows from the central limit theorem and the formulas $\mathbb{E}[Y_{r,i}]=0$ and $\mathrm{Var}(Y_{r,i})=\sigma^2+2\theta^2$ (see (\ref{mean}) and (\ref{var})). 

\vspace{2mm}

\noindent{4.} Let $(X,Y)$ be a bivariate normal random vector with  zero mean vector, variances $(\sigma_X^2,\sigma_Y^2)$ and correlation coefficient $\rho$. Denote the product of these correlated normal random variables by $Z=XY$. We will denote this distribution by $\mathrm{PN}(\rho,\sigma_X,\sigma_Y)$. Consider also the mean $\overline{Z}_n=n^{-1}(Z_1+\cdots+Z_n)$, where $Z_1,\ldots,Z_n$ are independent copies of $Z$. It was noted in the thesis of Gaunt \cite{gaunt thesis} that $Z\sim\mathrm{VG}(1,\rho \sigma_X\sigma_Y,\sigma_X\sigma_Y\sqrt{1-\rho^2},0)$ and later by Gaunt \cite{gaunt prod} that also 
\begin{equation}\label{vgrep}\overline{Z}_n\sim\mathrm{VG}(n,\rho \sigma_X\sigma_Y/n,\sigma_X\sigma_Y\sqrt{1-\rho^2}/n,0).
\end{equation}
To see this, suppose $\sigma_X=\sigma_Y=1$, with the general case following by rescaling. One can readily verify that $X$ and $W=(Y-\rho X)/\sqrt{1-\rho^2}$ are independent $N(0,1)$ random variables. Thus, $Z=XY=X(\sqrt{1-\rho^2}W+\rho X)=\sqrt{1-\rho^2}XW+\rho X^2$, so that $Z\sim \mathrm{VG}(1,\rho,\sqrt{1-\rho^2},0)$ by (\ref{gauntrep}). Finally, we deduce (\ref{vgrep}) from properties (\ref{muc}) and (\ref{conv}).

With (\ref{vgrep}), Gaunt \cite{gaunt prod} proved that the PDF of $\overline{Z}_n$ (previously given by independently by Mangilli, Plaszczynski and Tristram \cite{man15} and Nadarajah and Pog\'{a}ny \cite{np16}) is given by
\begin{equation*}p(x)=\frac{n^{(n+1)/2}2^{(1-n)/2}|x|^{(n-1)/2}}{s^{(n+1)/2}\sqrt{\pi(1-\rho^2)}\Gamma(n/2)}\exp\bigg(\frac{\rho n x}{s(1-\rho^2)} \bigg)K_{\frac{n-1}{2}}\bigg(\frac{n |x|}{s(1-\rho^2)}\bigg), \quad x\in\mathbb{R},
\end{equation*}
where $s=\sigma_X\sigma_Y$, with the PDF for $Z$ following on setting $n=1$. For an account of the distributional theory of $Z$ and $\overline{Z}$, see Gaunt \cite{gaunt22}. We note that the product of two normal random variables with non-zero means does not follow the VG distribution; the rather complicated exact formula for the PDF is given in Cui et al.\ \cite{cui}. 



\subsubsection{Superclasses}

\noindent{1.} The generalized hyperbolic (GH) distribution was introduced by Barndorff-Nielsen \cite{barndorff, barndorff2}, who studied it in the context of modelling dune movements.  Like the VG distribution, GH distributions are widely used in financial modelling; see for example, Bibby and S{\o}rensen \cite{bibby},  Eberlein and Keller \cite{eberlein1} and Eberlein and Prause \cite{eberlein2}. Properties of the GH distribution are given in Bibby and S{\o}rensen \cite{bibby} and Hammerstein \cite{hammersteinphd}. The PDF is
\begin{equation}\label{ghpdf}p(x)=\frac{(\gamma/\delta)^\lambda}{\sqrt{2\pi}K_\lambda(\delta\gamma)}\mathrm{e}^{\beta(x-\mu)}\frac{K_{\lambda-1/2}(\alpha\sqrt{\delta^2+(x-\mu)^2})}{(\sqrt{\delta^2+(x-\mu)^2}/\alpha)^{1/2-\lambda}}, \quad x\in\mathbb{R},
\end{equation}
where $\gamma=\sqrt{\alpha^2-\beta^2}$. 
The parameter domain is given by
\begin{align*}&\delta\geq0, \quad \gamma>0, \quad \lambda>0, \\
&\delta>0, \quad \gamma>0, \quad \lambda=0, \\
&\delta>0, \quad \gamma\geq0, \quad \lambda<0,
\end{align*}
and in each case $\mu\in\mathbb{R}$.  If $\delta=0$ or $\gamma=0$,  the PDF (\ref{ghpdf}) is defined as the limit obtained by using (\ref{Ktend0}).  In particular, taking $\delta\downarrow0$ in (\ref{ghpdf}) yields the density (\ref{vgpam3}) of the VG distribution in the parametrisation of Bibby and S{\o}rensen \cite{bibby}. A detailed study of limiting cases of GH distributions is given by Eberlein and Hammerstein \cite{eberlein}.

\vspace{2mm}

\noindent{2.} The VG distribution is a special case of the CGMY distribution of Carr et al.\ \cite{cgmy}, which was introduced in the context of financial modelling. The distribution is defined by its
 L\'{e}vy-Hinchin representation $\varphi(t)=\exp\big(
 \int_\mathbb{R}(\mathrm{e}^{\mathrm{i}tx}-1)\nu_{CGMY}(x)\,\mathrm{d}x \big),$
where the L\'evy density is given by
\begin{equation*}\nu_{CGMY}(x)=\frac{C}{|x|^{1+Y}}\mathrm{e}^{-G|x|}\mathbf{1}_{x<0}+\frac{C}{x^{1+Y}}\mathrm{e}^{-Mx}\mathbf{1}_{x>0}.
\end{equation*}
The VG distribution corresponds to the case $Y=0$. A closed-form formula is not available for the PDF of the CGMY distribution, although the characteristic function takes an elementary form (see Carr et al.\ \cite{cgmy}):
\begin{equation*}\varphi(t)=\exp\big(
C\Gamma(-Y)\big\{(M-\mathrm{i}t)^Y-M^Y+(G-\mathrm{i}t)^Y-G^Y\big\}\big),
\end{equation*}
from which one readily obtains formulas for lower order moments, as well the variance, skewness and kurtosis (again, see Carr et al.\ \cite{cgmy}).

Being a subclass of the CGMY distributions, the VG distributions form a subclass of the tempered stable distributions that was introduced by Koponen \cite{k95}. 
For properties and references for this class of distributions see, for example, K\"uchler and Tappe \cite{kt13}.

\vspace{2mm}

\noindent{3.} Consider the distribution
of the random variable $X = Z + L$, with independent $Z \sim N(0, \sigma^2)$ and $L$ is asymmetric Laplace. This distribution, introduced by Reed \cite{Reed1}, called the {\it normal-asymmetric Laplace}, has four parameters. Some basic distributional properties are given in Reed \cite{Reed1}, such as infinite divisibility, exact formulas for the PDF, CDF, moment and cumulants, and estimation methods are also discussed.

More generally, if $L$ is VG distributed (generalized asymmetric Laplace), then $X$ is said to be {\it generalized normal-Laplace} (GNL) distributed, introduced in Reed \cite{Reed2}. This distribution has five parameters. The article of Reed \cite{Reed2} includes financial applications. In the thesis of Wu \cite{MSthesis}, one can find a systematic treatment and bibliography for GNL distributions. 

\vspace{2mm} 

{\noindent{4.} 
A multivariate generalisation of the VG distribution in the symmetric case was first proposed by Madan and Seneta \cite{madan}. 
A number of basic properties of the (general non-symmetric) multivariate VG distribution are reviewed by Kozubowski, Podg\'{o}rski and Rychlik \cite{kpr13} (they refer to it as the multivariate generalized Laplace distribution). For an account of the matrix variate VG distribution, see Kozubowski, Mazur and Podg\'{o}rski \cite{kmp22} and references therein. Here, we collect just a few of the most fundamental properties of the multivariate VG distribution, noting that they generalise properties of the univariate VG distribution in a natural manner.

The multivariate VG distribution is often defined through its characteristic function. In line with Kozubowski et al.\ \cite{kpr13}, we say that a random vector $\mathbf{X}$ in $\mathbb{R}^d$ follows the multivariate VG distribution if its characteristic function is given by
\begin{equation}\label{mvvvcf}\varphi(\mathbf{t}) = \mathbb E\big[\mathrm{e}^{\mathrm{i}\mathbf{t}^\intercal \mathbf{X}}\big] =\mathrm{e}^{\mathrm{i}\boldsymbol{\mu}^\intercal \mathbf{t}}\big(1  - \mathrm{i}\boldsymbol{\theta}^\intercal\mathbf{t}+ \tfrac{1}{2}\mathbf{t}^\intercal\boldsymbol{\Sigma}\mathbf{t}\big)^{-s} , \quad \mathbf{t}\in\mathbb{R}^d,
\end{equation}
where $s>0$, $\boldsymbol{\theta}\in\mathbb{R}^d$, $\boldsymbol{\mu}\in\mathbb{R}^d$, and $\boldsymbol{\Sigma}$ is $d\times d$ non-negative definite symmetric matrix. Note that when $d=1$, the characteristic function (\ref{mvvvcf}) reduces to the characteristic function (\ref{cfcf}) of the univariate VG distribution (under a slightly different parametrisation). We can represent the multivariate VG distribution as a normal variance-mean mixture
\[
\mathbf{X} =_d \boldsymbol{\mu} + \boldsymbol{\theta}G + \sqrt{G}\mathbf{Z},
\]
for independent $G \sim \Gamma(s, 1)$ and $\mathbf{Z} \sim  N_d(\mathbf{0}, \boldsymbol{\Sigma})$ (the $d$-dimensional multivariate normal distribution with mean vector $\mathbf{0}$ and covariance matrix $\boldsymbol{\Sigma}$).

If the matrix $\boldsymbol{\Sigma}$ is positive-definite, then the distribution is non-degenerate, that is truly $d$-dimensional, with PDF
\[
p(\mathbf{x}) = \frac{2\exp(\boldsymbol{\theta}^\intercal\boldsymbol{\Sigma}^{-1}(\mathbf{x}-\boldsymbol{\mu}))}{(2\pi)^{d/2}\Gamma(s)|\boldsymbol{\Sigma}|^{1/2}}\left(\frac{Q(\mathbf{x},\boldsymbol{\mu})}{C(\boldsymbol{\Sigma}, \boldsymbol{\theta})}\right)^{s-d/2}K_{s-d/2}(Q(\mathbf{x},\boldsymbol{\mu})C(\boldsymbol{\Sigma},\boldsymbol{\theta})), \quad \mathbf{x}\in\mathbb{R}^d,\]
where $Q(\mathbf{x},\boldsymbol{\mu})=((\mathbf{x}-\boldsymbol{\mu})^\intercal\boldsymbol{\Sigma}^{-1}(\mathbf{x}-\boldsymbol{\mu}))^{1/2}$ and $C(\boldsymbol{\Sigma},\boldsymbol{\theta})=(2+\boldsymbol{\theta}^\intercal\boldsymbol{\Sigma}^{-1}\boldsymbol{\theta})^{1/2}$.

\subsubsection{Products and quotients of VG random variables}

Let $(X_i)_{1\leq i\leq n}\sim\mathrm{VG}(r_i,0,\sigma_i,0)$ be independent, and set $Z_n=\prod_{i=1}^nX_i$ and $\sigma=\sigma_1\cdots\sigma_n$. A formula for the PDF of the product $Z_n$ was obtained by Gaunt, Mijoule and Swan \cite{gms}:
\begin{equation*}
f_{Z_n}(x)=\frac{1}{2^n\pi^{n/2}\sigma}\prod_{j=1}^n\frac{1}{\Gamma(r_j/2)}G^{2n,0}_{0,2n}\bigg(\frac{x^2}{2^{2n}\sigma^2}\, \bigg| \, \frac{r_1-1}{2},\ldots,\frac{r_n-1}{2},0,\ldots,0 \bigg), \quad x\in\mathbb{R},
\end{equation*}
where $G_{0,k}^{k,0}(x\,|\,)$ is a Meijer $G$-function (see Chapter 16 of Olver et al.\ \cite{olver} for a definition and basic properties). We are not aware of an exact formula for the product of two or more independent $\mathrm{VG}(r,\theta,\sigma,0)$ random variables with $\theta\not=0$.

Suppose that $X\sim \mathrm{VG}(r,\theta_1,\sigma_1,0)$ and $Y\sim\mathrm{VG}(s,\theta_2,\sigma_2,0)$ are independent, where $r,s>0$, $\theta_i\in\mathbb{R}$, $\sigma_i>0$, $i=1,2$. An exact formula for the PDF of ratio $Z=X/Y$ in terms of an infinite series involving Gaussian hypergeometric functions (See Chapter 15 of Olver et al.\ \cite{olver} for a definition and basic properties) was obtained by Gaunt and Li \cite{gl23}. Previously, the following simpler exact formula, expressed as a single Gaussian hypergeometric function, for the case $\theta_1=\theta_2=0$ and $\sigma_1=\sigma_2=1$  was obtained by Nadarajah and Kotz \cite{nk06}.
For $x\in\mathbb{R}$,
\begin{align*}f_Z(x)=\frac{2|x|^{-s-1}\Gamma((r+1)/2)\Gamma((s+1)/2)}{\pi(r+s)\Gamma(r/2)\Gamma(s/2)}{}_{2}F_1\bigg(\frac{r+s}{2},\frac{s+1}{2};\frac{r+s}{2}+1;1-x^{-2}\bigg).
\end{align*}
We refer the reader to the two aforementioned works for further properties of the VG ratio distribution, such as formulas for the PDF expressed in terms of elementary functions for particular parameter values, a formula for the CDF, fractional moments, and asymptotic behaviour of the PDF and tail probabilities, from which it is seen that the mean does not exist. 

 
\subsection{Stein characterisations} \label{sec2.6}

The following Stein characterisation of the VG distribution was given by Gaunt \cite{gaunt vg}. Let $W$ be a real-valued random variable. Then $W\sim\mathrm{VG}(r,\theta,\sigma,\mu)$ if and only if
\begin{align}\label{char1}\mathbb{E}\big[\sigma^2(W-\mu)g''(W)+(\sigma^2r+2\theta(W-\mu))g'(W)
+(r\theta-(W-\mu))g(W)\big]&=0
\end{align}
for all twice differentiable $g:\mathbb{R}\rightarrow\mathbb{R}$ for which the expectations $\mathbb{E}|g(X)|$, $\mathbb{E}|Xg(X)|$, $\mathbb{E}|g'(X)|$, $\mathbb{E}|Xg'(X)|$ and $\mathbb{E}|Xg''(X)|$ are finite for $X\sim\mathrm{VG}(r,\theta,\sigma,\mu)$. 

Here we provide a sketch of the argument. We start with necessity. The PDF (\ref{vgdefn}) of the $\mathrm{VG}(r,\theta,\sigma,\mu)$ distribution satisfies the ordinary differential equation (ODE)
\begin{equation}\label{pdfode}\mathcal{L}p(x):=\sigma^2xp''(x)-(\sigma^2(r-2)+2\theta x)p'(x)+(\theta(r-2)-x)p(x)=0,
\end{equation}
which can be read off from a more general ODE satisfied by the PDF of the GH distribution given in Gaunt \cite[Corollary 3.3]{gaunt gh}. From (\ref{pdfode}) we have that $\int_{-\infty}^\infty g(x)\mathcal{L}p(x)\,\mathrm{d}x=0$ and integrating by parts twice yields (\ref{char1}).
 To prove sufficiency, suppose $W$ is a real-valued random variable. For ease of exposition, we set $\mu=0$, with the general case following by a simple translation. Taking $g(x)=\mathrm{e}^{\mathrm{i}tx}$
 (a twice differentiable and bounded function) 
 in (\ref{char1}) and letting $\varphi(t)=\mathbb{E}[\mathrm{e}^{\mathrm{i}tW}]$ leads to the ODE
\begin{equation}\label{odee}(\sigma^2t^2-2\mathrm{i}\theta t+1)\varphi'(t)+(\mathrm{i}\sigma^2rt+r\theta)\varphi(t)=0.
\end{equation}
Here we applied the characterising equation to the real and imaginary parts of $g(x) = \mathrm{e}^{\mathrm{i}tx}$. Solving (\ref{odee}) subject to the condition $\varphi(0)=1$ gives that $\varphi(t)=(1-2\mathrm{i}\theta t+\sigma^2t^2)^{-r/2}$, in agreement with the characteristic function (\ref{cfcf}) when $\mu=0$, and hence $W\sim\mathrm{VG}(r,\theta,\sigma,0)$.

Stein characterisations are most commonly used in Stein's method (Stein \cite{stein}) to bound the distance between two probability distributions with respect to a probability metric; however, one can also use them to establish distributional properties. Indeed, letting $Y\sim\mathrm{VG}(r,\theta,\sigma,0)$ and substituting $g_1(x)=x^k$ and $g_2(x)=(x-\mathbb{E}[Y])^k=(x-r\theta)^k$ into (\ref{char1}) (with $\mu=0$) leads to the following recursions for the $k$-th raw moment $\mu_k'=\mathbb{E}[Y^k]$ and the $k$-th central moment $\mu_k=\mathbb{E}[(Y-\mathbb{E}[Y])^k]$:
\begin{align}\label{rec1}\mu_{k+1}'&=\theta(2k+r)\mu_k'+\sigma^2k(r+k-1)\mu_{k-1}', \quad k\geq1, \\
\label{rec2}\mu_{k+1}&=2k\theta\mu_k+k\big(\sigma^2(r+k-1)+2\theta^2r\big)\mu_{k-1}+k(k-1)r\theta\sigma^2\mu_{k-2}, \quad k\geq2.
\end{align}
Lower order raw (central) moments can be efficiently computed using the recurrences given just the first raw (first and second central) moments, along with the basic condition $\mu_0=1$.

\subsection{Moments and cumulants}\label{sec2.7}

The mean and variance of $X\sim\mathrm{VG}(r,\theta,\sigma,\mu)$ are readily calculated
 via the standard method of obtaining lower order moments from moment generating functions (using (\ref{mgf})).
We have
\begin{align}\label{mean}\mathbb{E}[X]&=\mu+r\theta, \\
\label{var}\mathrm{Var}(X)&=r(\sigma^2+2\theta^2).
\end{align}
Now, let $Y\sim\mathrm{VG}(r,\theta,\sigma,0)$. Lower order raw and central moments of $Y$ can be efficiently obtained from the recurrences (\ref{rec1}) and (\ref{rec2}). The central moments of $X\sim\mathrm{VG}(r,\theta,\sigma,\mu)$ are equal to the central moments of $Y$, and raw moments for $X$ are readily obtained through the formula $\mathbb{E}[X^k]=\sum_{j=0}^k\binom{k}{j}\mu^j\mathbb{E}[Y^k]$. The first four raw moments of $Y\sim\mathrm{VG}(r,\theta,\sigma,0)$ are
\begin{align*}\mu_1'&=r\theta, \quad
\mu_2'=r(\sigma^2+(r+2)\theta^2), \quad
\mu_3'=r(r+2)\theta\big(3\sigma^2+(r+4)\theta^2\big), \\
\mu_4'&=r(r+2)\big(3\sigma^4+6(r+4)\theta^2\sigma^2+(r+4)(r+6)\theta^4\big),
\end{align*}
whilst the first four central moments are
\begin{align}\label{mom1}\mu_1&=0,\\
\label{mom2}\mu_2&=r(\sigma^2+2\theta^2),  \\
\label{mom3}\mu_3&=2r\theta(3\sigma^2+4\theta^2), \\
\label{mom4}\mu_4&=3r\big((r+2)\sigma^4+(4r+16)\theta^2\sigma^2+(4r+16)\theta^4\big).
\end{align}
The skewness $\gamma_1=\mu_3/\mu_2^{3/2}$ and kurtosis $\beta_2=\mu_4/\mu_2^2$ are thus given by
\begin{align*}\gamma_1=\frac{2\theta(3\sigma^2+4\theta^2)}{\sqrt{r}(\sigma^2+2\theta^2)^{3/2}}, \quad
\beta_2=\frac{3((r+2)\sigma^4+(4r+16)\theta^2\sigma^2+(4r+16)\theta^4)}{r(\sigma^2+2\theta^2)^2},
\end{align*}
and the excess kurtosis $\gamma_2=\mu_4/\mu_2^2-3$ is equal to
\[\gamma_2=\frac{6(\sigma^4+8\theta^2\sigma^2+8\theta^4)}{r(\sigma^2+2\theta^2)^2}.\]
Formulas for the first four central moments of the VG distribution, as well as the skewness and kurtosis, were given by Seneta \cite{s04}.


Higher order moments can be calculated using the representation (\ref{gamrep}) (with $\mu=0$) of the $\mathrm{VG}(r,\theta,\sigma,0)$ distribution and the formula $\mathbb{E}[V^k]=\lambda^{-k}\Gamma(r/2+k)/\Gamma(r/2)
$, $k\geq1$, where $V\sim\Gamma(r/2,\lambda)$. For $k\geq1$,
\begin{align*}\mu_k'
=\frac{\sigma^{2k}}{(\Gamma(r/2))^2}\sum_{j=0}^k\binom{k}{j}(-\lambda_-)^j\lambda_+^{k-j}\Gamma\Big(\frac{r}{2}+j\Big)\Gamma\Big(\frac{r}{2}+k-j\Big)
.
\end{align*}
We now return to the general case $r>0$, for which the following expressions, involving the hypergeometric function, for the raw and absolute moments of $Y\sim\mathrm{VG}(r,\theta,\sigma,0)$ were obtained by Gaunt \cite{gauntvgmom}. Let $k\in\mathbb{Z}^+$, and define $\ell:=\lceil k/2\rceil+1/2$ and $m:=k\,\mathrm{mod}\,2$. Then
\begin{equation*}\mu_k'=\frac{2^{k+m}\theta^m\sigma^{r+2k}}{\sqrt{\pi}(\theta^2+\sigma^2)^{(r+k+m)/2}\Gamma(r/2)}\Gamma\Big(\frac{r-1}{2}+\ell\Big)\Gamma(\ell)\,{}_2F_1\bigg(\ell,\frac{r-1}{2}+\ell;\frac{1}{2}+m;\frac{\theta^2}{\theta^2+\sigma^2}\bigg),
\end{equation*}
whilst, for $k>k_*=k_*(r):=\mathrm{max}\{-1,-r\}$,
\begin{align*}\mathbb{E}[|Y|^k]=\frac{2^k\sigma^{r+2k}}{\sqrt{\pi}(\theta^2+\sigma^2)^{(r+k)/2}\Gamma(r/2)}\Gamma\Big(\frac{r+k}{2}\Big)\Gamma\Big(\frac{k+1}{2}\Big)\,{}_2F_1\bigg(\frac{k+1}{2},\frac{r+k}{2};\frac{1}{2};\frac{\theta^2}{\theta^2+\sigma^2}\bigg).
\end{align*}
As ${}_2F_1(a,b;c;0)=1$, we observe that in the case $\theta=0$ the raw moments simplify to
\begin{equation*}\mu_k'=\begin{cases} \displaystyle \frac{2^{k}\sigma^k}{\sqrt{\pi}}\frac{\Gamma((r+k)/2)\Gamma((k+1)/2)}{\Gamma(r/2)}, & \quad \text{if $k$ is even}, \\
\displaystyle 0, & \quad \text{if $k$ is odd}. \end{cases}
\end{equation*}

Central moments of general order of the $\mathrm{VG}(r,\theta,\sigma,0)$ distribution can be calculated using the representation (\ref{gamrep}) and the formula $\mathbb{E}[(V-\mathbb{E}[V])^k]=\lambda^{-k}U(-k,1-k-r/2,-r/2)$, where $V\sim\Gamma(r/2,\lambda)$ and $U(a,b,x)$ is a confluent hypergeometric function of the second kind. This formula is given for the case that $V\sim \chi_r^2$ by Weisstein \cite{mathworld}, and the generalisation to the gamma distribution is obvious. We observe that $U(a,b,x)$ is a polynomial when $a=-m$, $m=0,1,2,\ldots$:
$U(-m,b,x)=(-1)^m\sum_{j=0}^m\binom{m}{j}(b+j)_{m-j}(-x)^j$. Here the Pochhammer symbol is given by $(a)_0=0$ and $(a)_j=a(a+1)(a+2)\cdots(a+j-1)$, $j\geq1$ (see Olver et al.\ \cite[Section 13.2(i)]{olver}). For $k\geq1$,
\begin{align*}\mu_k
=\sigma^{2k}\sum_{j=0}^k\binom{k}{j}(-\lambda_-)^j\lambda_+^{k-j}U\bigg(-j,1-j-\frac{r}{2},-\frac{r}{2}\bigg)U\bigg(j-k,1+j-k-\frac{r}{2},-\frac{r}{2}\bigg).
\end{align*}   
Other formulas for the 
raw and central moments of the VG distribution are given by Scott \cite{scott}.

 The cumulants of the $\mathrm{VG}(r,\theta,\sigma,\mu)$ distribution can be derived by Taylor expanding the logarithms in (\ref{kfkf}) or using the representation (\ref{gamrep}). The latter approach involves using the facts that, for independent random variables $S_1$ and $S_2$ and constant $c$, the $k$-th cumulant enjoys the properties $\kappa_k(cS_1)=c^k\kappa_k(S_1)$ and $\kappa_k(S_1+S_2)=\kappa_k(S_1)+\kappa_k(S_2)$, and $\kappa_k(V)=\lambda^{-k}(k-1)!\,r$ for $V\sim\Gamma(r/2,\lambda)$. For $k\geq1$, 
\begin{align*}\kappa_k=\mu\delta_{k,1}+\frac{(k-1)!\,r}{2}\big[(\theta+\sqrt{\theta^2+\sigma^2})^k+(\theta-\sqrt{\theta^2+\sigma^2})^k\big],
\end{align*}
where $\delta_{i,j}$ is the Kronecker delta (see McKay \cite{m32}). In particular,
\begin{align*}\kappa_1&=\mu+r\theta, \quad\kappa_2=r(\sigma^2+2\theta^2), \quad\kappa_3=2r\theta (3\sigma^2+4\theta^2),\quad \kappa_4=6r(\sigma^4+8\theta^2\sigma^2+8\theta^4),\\
\kappa_5&=24r\theta(5\sigma^4+20\theta^2\sigma^2+16\theta^4), \quad \kappa_6=120r(\sigma^6+18\theta^2\sigma^4+48\theta^4\sigma^2+32
\theta^6).
\end{align*}

\subsection{Mode and median}\label{sec2.8}

A detailed study of the mode and median of the GH and VG distributions was undertaken by Gaunt and Merkle \cite{gm21}; here we provide a summary for the VG distribution. 

The $\mathrm{VG}(r,\theta,\sigma,\mu)$ distribution is unimodal. This follows since the distribution is self-decomposable and self-decomposable distributions are unimodal (Yamazato \cite{y78}). Denote the mode of the $\mathrm{VG}(r,\theta,\sigma,\mu)$ distribution by $M_r$. 
For $0<r\leq2$, $\theta\in\mathbb{R}$, $\sigma>0$, or $r>0$, $\theta=0$, $\sigma>0$, we have $M_r=\mu$. This is clear when $\theta=0$, because the $VG(r,0,\sigma,\mu)$ distribution is symmetric about $\mu$; for $0<r\leq 2$ and $\theta\in\mathbb{R}$ this follows because, for $0<\nu\leq1/2$ and $|\beta|<1$, the function $x\mapsto \mathrm{e}^{\beta x}|x|^\nu K_\nu(|x|)$ is increasing on $(-\infty,0)$ and decreasing on $(0,\infty)$ (this is easily inferred from part (i) of Lemma 5.1 of Gaunt \cite{gaunt vg3}). 

Suppose now that $r>2$, $\theta\in\mathbb{R}$, $\sigma>0$, $\mu\in\mathbb{R}$. 
 Using (\ref{ddbk}), we obtain that $M_r=\mu+\mathrm{sgn}(\theta)\cdot x^*$, where $x^*$ is the unique positive solution of the equation 
\begin{equation}\label{xstar}K_{\frac{r-3}{2}}\bigg(\frac{\sqrt{\theta^2+\sigma^2}}{\sigma^2}x\bigg)=\frac{|\theta|}{\sqrt{\theta^2+\sigma^2}}K_{\frac{r-1}{2}}\bigg(\frac{\sqrt{\theta^2+\sigma^2}}{\sigma^2}x\bigg).
\end{equation}
When $r=4$ and $r=6$, applying (\ref{special}) to (\ref{xstar}) leads to simple algebraic equations for $x^*$ which when solved give
\[M_4=\mu+\theta\bigg(1+\frac{1}{\sqrt{1+\kappa}}\bigg), \quad
M_6 =\mu+\frac{\theta}{2}\bigg(1+\frac{1}{\sqrt{1+\kappa}}\bigg)\bigg(3-\sqrt{1+\kappa}+\sqrt{6\sqrt{1+\kappa}+\kappa-2}\bigg),\]
where $\kappa=\sigma^2/\theta^2$.
For $r=8$ and $r=10$, we can apply (\ref{special}) to (\ref{xstar}) to obtain cubic and quartic equations for $x^*$; however, their solutions $M_8$ and $M_{10}$ are too complicated to be worth reporting. For other values of $r$, exact formulas for $M_r$ are not available, except for the case $\theta=0$ in which we have already noted that the mode is equal to $\mu$.


Exact formulas for $M_r$ are not available for general $r>2$; however, simple and accurate lower and upper bounds can be obtained.
Fix $\theta>0$; bounds for $\theta<0$ follow immediately as if $X\sim \mathrm{VG}(r,\theta,\sigma,\mu)$ then $-X\sim\mathrm{VG}(r,-\theta,\sigma,-\mu)$ (see (\ref{muc})). Applying the lower and upper bounds of (\ref{seg1}) to (\ref{xstar}) yields simple algebraic inequalities for $x^*$. Solving these inequalities and using that $M_r>\mu$, for $r>2$ and $\theta>0$, leads to the following double inequality:
\begin{equation}\label{mode1} \theta(r-3)_+<M_r-\mu<\theta(r-2), \quad r>2, 
\end{equation}
where $x_+=\max\{0,x\}$. Note that there is equality in the upper bound when $r=2$. Applying inequality (\ref{seg3}) to equation (\ref{xstar}) gives that 
\begin{equation}\label{mode2}M_r>\mu+\frac{\theta}{2}\bigg[r-2+\sqrt{\frac{\theta^2(r-2)^2+\sigma^2(r-4)^2}{\theta^2+\sigma^2}}\bigg], \quad r>4,
\end{equation}
with equality if and only if $r=4$, and the inequality is reversed if $2<r<4$. 
Inequality (\ref{mode2}) improves on the lower bound in (\ref{mode1}) in its range of validity $r>4$.  The reversed inequality (\ref{mode2}) is more accurate than the upper bound in (\ref{mode1}) for $3<r<4$, but the reverse is true for $2<r<3$. 


Let $X\sim\mathrm{VG}(r,\theta,\sigma,\mu)$. We note that exact formulas for the median are available for certain parameter values: $\mathrm{Med}(X)=\mu$ if $\theta=0$, whilst, for $r=2$,
\begin{equation}\label{medform}\mathrm{Med}(X)=\mu+\mathrm{sgn}(\theta)\cdot\big(|\theta|+\sqrt{\theta^2+\sigma^2}\big)\log\bigg(1+\frac{|\theta|}{\sqrt{\theta^2+\sigma^2}}\bigg),
\end{equation}
which follows from the formula for the  median of the asymmetric Laplace distribution (see Kozubowski and Podg\'{o}rski \cite{kp01}). For other parameter values, an exact closed-form formula is not available. Moreover, in contrast to the mode, accurate lower and upper bounds for the median are not available in the literature. 
Gaunt and Merkle \cite{gm21} have conjectured accurate bounds for the median, which we present below.  The numerical results of Table \ref{table1}, which are taken from Gaunt and Merkle \cite{gm21}, support the conjectured bounds, since each entry in the table lies between the conjectured lower and upper bounds. 

\begin{conjecture}\label{conj3}Let $X\sim\mathrm{VG}(r,\theta,\sigma,\mu)$ with $r,\theta,\sigma>0$. Then it is conjectured that
\begin{equation*}\mu+(r-1)\theta<\mathrm{Med}(X)<\mu+r\theta \mathrm{e}^{-2/3r}<\mu+\bigg(r-\frac{2}{3}+\frac{2}{9r}\bigg)\theta, \quad r>0,
\end{equation*}
and 
\begin{equation*}\mathrm{Med}(X)\leq\mu+(r+2\log2-2)\theta, \quad r\geq2.
\end{equation*}
\end{conjecture}

\begin{table}[h]
\centering
\caption{\footnotesize{Median of the $\mathrm{VG}(r,1,\sigma,0)$ distribution.}}
\label{table1}
{\normalsize
\begin{tabular}{|c|rrrrrr|}
\hline
 \backslashbox{$r$}{$\sigma$}       &    0.1 &    0.3 &    1 & 3 &    10 &    30   \\
 \hline
0.5  & 0.0863 & 0.0798 & 0.0502 & 0.0195 & 0.00582 & 0.00192 \\
1  & 0.454 & 0.444 & 0.380 & 0.276 & 0.198  &  0.157 \\
2.5  & 1.872 & 1.861 & 1.775 & 1.621 & 1.531 & 1.507  \\
5 & 4.350 & 4.338 & 4.246 & 4.084 & 4.012 & 4.001  \\ 
10 & 9.340 & 9.328  & 9.233 & 9.071 & 9.009 & 9.001  \\ 
  \hline
\end{tabular}}
\end{table}

Since $\mathbb{E}[X_r]=\mu+r\theta$ for $X_r\sim\mathrm{VG}(r,\theta,\sigma,\mu)$, it follows from  (\ref{mode1}) and the fact that $M_r=\mu$ for $0<r\leq2$ that, if $\theta>0$, then $M_r\leq\mathbb{E}[X_r]$ for all $r>0$. 
Moreover, if the conjectured median bounds hold, then it would follow that $M_r\leq\mathrm{Med}(X_r)\leq\mathbb{E}[X_r]$ for all $r,\theta,\sigma>0$, meaning that the $\mathrm{VG}(r,\theta,\sigma,\mu)$ distribution would satisfy the mean-median-mode inequality (Groeneveld and Meeden \cite{gm77} and van Zwet \cite{v79}). 
 
\section{Parameter Estimation}\label{sec3}
Due to its numerous applications, parameter estimation of the VG distribution is of much interest. In this section, we give an overview of 
available techniques from the literature.
In the following, we let ${\bf X}=(X_1,\ldots,X_n)$ be an i.i.d.\ sample of $\mathrm{VG}(r, \theta,\sigma,\mu)$ random variables from which we wish to estimate the four unknown parameters $r, \theta,\sigma$ and $\mu$.

\subsection{Method of moments estimation}
Moment estimation consists of estimating the raw or central moments through the arithmetic means, i.e. $\hat{\mu}_j'=n^{-1}\sum_{i=1}^n X_i^j$, $j=1,\ldots,4$, and $\hat{\mu}_j=n^{-1}\sum_{i=1}^n (X_i- \hat{\mu}_j')^j$, $j=1,\ldots,4$.
As noted by Finlay and Seneta \cite{fs08}, in the four parameter case, it is not possible to resolve $\hat{\mu}_j=\mu_j$, $j=1,\ldots,4$, for $(r, \theta,\sigma,\mu)$ explicitly. However, in th symmetric case $\theta=0$, this is possible. To illustrate this, we will work in the parameterisation (\ref{vgpam2}) of Madan and Seneta \cite{madan}.
Using the formulas for the first, second and fourth central moments of the VG distribution (\ref{mom1}), (\ref{mom2}) and (\ref{mom4}), one readily obtains the following explicit consistent estimators, which is also unbiased in the case of the estimator for $\mu$ (Madan and Seneta \cite{madan}):
\begin{align*}\hat{\mu}=\frac{1}{n}\sum_{i=1}^nX_i, \quad \hat{\sigma}_0^2=\frac{1}{n}\sum_{i=1}^n(X_i-\overline{X})^2, \quad \hat{\nu}=\frac{1}{3\hat{\sigma}_0^4 n}\bigg\{\sum_{i=1}^n(X_i-\overline{X})^4\bigg\}-1.
\end{align*}
 In Finlay and Seneta \cite{fs08}, the moment estimates are obtained by minimising the quantity $ \sum_{j=1}^4 (( \hat{\mu}_j'- \mu_j')/ \hat{\mu}_j')^2$
with respect to $(r, \theta,\sigma,\mu)$, where $\mu_j'$, $j=1,\ldots,4$, are as in Section \ref{sec2.7}. In Seneta \cite{s04}, explicit moment estimators are obtained for the parametrisation \eqref{vgpam2} by assuming that expressions of the form $\theta_0^k$, $ k\geq 2$ are negligible since the true value of $\theta_0$ is often small in finanical modelling applications. An advantage of the moment estimators is the easy numerical implementation. A competitive simulation study in order to evaluate the performance of the moment estimators is carried out in Finlay and Seneta \cite{fs08}.

\subsection{Maximum-likelihood estimation}
We maximise the log-likelihood function
\begin{align*}
    \ell({\bf X},r, \theta,\sigma,\mu)= \sum_{i=1}^n \log p(X_i,r, \theta,\sigma,\mu)
\end{align*}
with respect to $(r, \theta,\sigma,\mu)$, where $p$ is the VG density (\ref{vgdefn}). Due to the complexity of the density involving the modified Bessel function of the second kind, there is no analytic solution to this optimisation problem and one relies on numerical methods in order to determine the global maximum. We also stress the density function admits a singularity at $\mu$ for $r\leq1$. There exist implementations of numerical procedures for maximum-likelihood estimation, such as the \textit{VarianceGamma} R package of Scott and Dong \cite{scott18} through the \textit{VGFit}-function. Madan and Seneta \cite{madan}  proposed a transformed maximum likelihood method which permits to perform optimisation with respect to a more tractable density function (see also Madan and Seneta \cite{ms87b}). In Madan et al.\ \cite{mcc98}, MLE is used in order to estimate the parameters for a real financial data example. In Bee, Dickson and Santi \cite{bee18}, the \textit{VarianceGamma} package is compared to the \textit{ghyp} R package of Breymann and L\"uthi \cite{breymann13}, which implements the more general GH distribution, and to the ECM algorithm that we describe in Section \ref{sec3.3}. In Cervellera and Tucci \cite{cer17}, the authors investigate three different numerical optimisation procedures (including the \textit{VarianceGamma} R package, \textit{Matlab} and \textit{Ezgrad}) and point out several difficulties that can be encountered when using MLE. More precisely, the authors exhibit that the likelihood function is sensitive to the removal of a single observation and that in some cases the algorithms produce different output for different starting values. Due to the non-differentiability of the likelihood function with respect to $\mu$, the variance-gamma distribution does not fulfil the classical regularity assumptions that are required to apply the standard asymptotic theory of maximum likelihood estimators. However, in special cases in which certain parameters are known, the asymptotic theory has already been studied (see, for example, Kotz et al.\ \cite{kkp01} for the Laplace and asymmetric Laplace distributions).

\subsection{ECM algorithms}\label{sec3.3}
Motivated by the complexity of the optimisation of the likelihood function, several approaches using ECM (Expectation-Conditional Maximisation) algorithms have been studied recently. The algorithm we present here was first introduced by Nitithumbundit and Chan \cite{niti15} and McNicholas, McNicholas and Browne \cite{mcni13} for the multivariate VG model and mixtures of the latter.  We restrict ourselves to the one-dimensional case and present the algorithm described in detail in Bee et al.\ \cite{bee18}. Here we will consider the VG model under the parametrisation \eqref{vgpam2}. ECM algorithms are considered as iterative procedures based on the maximisation of the likelihood function with missing data. Letting $S\sim\Gamma(\alpha,\alpha)$, 
$T\sim N(0,1)$ be independent and using \eqref{rep2} we know that 
\begin{align*}
X=\mu + \theta_0 S +\sigma_0 \sqrt{S} T \sim \mathrm{VG}_2(\alpha,\theta_0,\sigma_0,\mu).
\end{align*}
We now treat ${\bf S}=(S_1,\ldots ,S_n)$ as missing data because the conditional distribution of $X$ on $S$ is $N(\mu+\theta_0 S, \sigma_0^2S)$, which means we can write the respective joint log-likelihood function as a product of a conditional normal and a gamma density. More precisely,
\begin{align*}
\ell_c({\bf X},{\bf S},\alpha, \theta_0,\sigma_0,\mu)=&-\frac{n}{2}\log \big(\sqrt{r} \sigma_0^2\big)
-\frac{1}{2} \sum_{i=1}^n \frac{1}{S_i} \frac{(X_i-\mu - S_i \theta_0)^2}{\sigma_0^2} \\
& + n \alpha \log (\alpha) - n \log (\gamma(\alpha)) + (\alpha -1)\sum_{i=1}^n \log S_i - \alpha  \sum_{i=1}^n S_i \\
=& \ell_c^{(1)}({\bf X},{\bf S}, \theta_0,\sigma_0,\mu) + \ell_c^{(2)}({\bf X},{\bf S},\alpha).
\end{align*}
ECM-type algorithms iterate an expectation (E)-step and a conditional maximisation (CM)-step. 

\vspace{0.2cm}

\noindent{{\bf E-step:}} Here one computes the conditional expectation $\mathbb{E}[\ell_c({\bf X},{\bf S},\alpha, \theta_0,\sigma_0,\mu) \, \vert \, {\bf X} ]$, which breaks down to the computation of $\mathbb{E}[S_i \, \vert \, {\bf X} ]$, $\mathbb{E}[1/S_i \, \vert \, {\bf X} ]$ and $\mathbb{E}[\log S_i \, \vert \, {\bf X} ]$. Embrechts \cite{embr83} showed that the conditional distribution of $S$ on $X$ is generalized inverse Gaussian and hence we can estimate $S_i$, $1/S_i$ and $\log S_i$ by the corresponding conditional expectations
\begin{align}
    \hat{S_i}&=\mathbb{E}[S_i \, \vert \, {\bf X} ]= \frac{\delta_i K_{\alpha+1/2}\big(\sqrt{2\alpha+\theta_0^2/\sigma_0^2}\delta_i \big)}{\sqrt{2\alpha+\theta_0^2/\sigma_0^2}K_{\alpha-1/2}\big(\sqrt{2\alpha+\theta_0^2/\sigma_0^2}\delta_i \big)}, \label{em_s1} \\
     \widehat{1/S_i}&=\mathbb{E}[1/S_i \, \vert \, {\bf X} ]= \frac{\sqrt{2\alpha+\theta_0^2/\sigma_0^2} K_{\alpha-3/2}\big(\sqrt{2\alpha+\theta_0^2/\sigma_0^2}\delta_i \big)}{\delta_iK_{\alpha-1/2}\big(\sqrt{2\alpha+\theta_0^2/\sigma_0^2}\delta_i \big)} \label{em_s2}, \\
      \widehat{\log S_i}&=\mathbb{E}[\log S_i \, \vert \, {\bf X} ]= \log \bigg(\frac{\delta_i}{\sqrt{2\alpha+\theta_0^2/\sigma_0^2}} \bigg)+ \frac{ \hat{K}_{\alpha-1/2}^{(1,0)}\big(\sqrt{2\alpha+\theta_0^2/\sigma_0^2}\delta_i \big)}{K_{\alpha-1/2}\big(\sqrt{2\alpha+\theta_0^2/\sigma_0^2}\delta_i \big)}, \label{em_s3}
\end{align}
where  $\delta_i=(X_i-\mu)^2/\sigma_0^2$ and $\hat{K}_{\alpha}^{(1,0)}(z)=(K_{\alpha+h}(z)-K_{\alpha-h}(z))/(2h)$, for some small constant $h>0$, is an approximation of $\partial K_{\tau}(z) / \partial \tau \vert_{\tau=\alpha}$ (see Bee et al.\ \cite[p.\ 76]{bee18}). 

\vspace{0.2cm} 

\noindent{{\bf CM-step:}} In this step, one maximises the conditional expectation of the log-likelihood function $\mathbb{E}[\ell_c({\bf X},{\bf S},\alpha, \theta_0,\sigma_0,\mu)\, \vert \ {\bf X }]$ with respect to the parameter vector $(\alpha,\theta_0,\sigma_0,\mu)$. One can maximise $\mathbb{E}[\ell_c^{(1)}({\bf X},{\bf S},\theta_0,\sigma_0,\mu)\, \vert \ {\bf X }]$ and $\mathbb{E}[\ell_c^{(2)}({\bf X},{\bf S},\alpha)\, \vert \ {\bf X }]$ separately. The first optimisation problem can be solved analytically by
\begin{align}
    \hat{\mu}&=\frac{\sum X_i \widehat{1/S_i} \sum \hat{S}_i - n \sum X_i}{\sum \widehat{1/S_i} \sum \hat{S}_i -n^2}, \label{emparam1} \\
    \hat{\theta}_0&=\frac{\sum X_i  - n \hat{\mu}}{\hat{S}_i}, \label{emparam2}  \\
    \hat{\sigma}^2_0&=\frac{1}{n}\sum \widehat{1/S_i} (X_i-\hat{\mu})^2-\frac{1}{n}\hat{\theta}_0^2 \sum \hat{S}_i. \label{emparam3} 
\end{align}
For optimisation with respect to $\alpha$, no analytic solution is available. It is suggested to optimise either $\mathbb{E}[\ell_c^{(2)}({\bf X},{\bf S},\alpha)\, \vert \ {\bf X }]$ or the likelihood of the VG distribution with updated parameters $\ell({\bf X},{\bf S},\alpha, \hat{\theta}_0,\hat{\sigma}_0,\hat{\mu})$. Bee et al.\ \cite{bee18} carried out a competitive simulation study, finding the latter of the two options to be preferable. Given a vector of starting values $(\hat{\alpha}^{(0)},\hat{\theta}_0^{(0)},\hat{\sigma}_0^{(0)},\hat{\mu}^{(0)})$, a pseudo-code for the $t$-th iteration of the algorithm is:
\begin{itemize}
    \item {\bf E-step 1:} Compute $\widehat{S_i}^{(t-1/2)}$ and $\widehat{1/S_i}^{(t)}$ using \eqref{em_s1} and \eqref{em_s2} and the current estimates $\hat{\alpha}^{(t-1)}$, $\hat{\theta}_0^{(t-1)}$, $\hat{\sigma}_0^{(t-1)}$ and $\hat{\mu}^{(t-1)}$.
    \item {\bf CM-step 1:} Compute $\hat{\theta}_0^{(t)}$, $\hat{\sigma}_0^{(t)}$ and $\hat{\mu}^{(t)}$ using \eqref{emparam1}, \eqref{emparam2} and \eqref{emparam3}.
    \item {\bf E-step 2:} Compute $\widehat{S_i}^{(t)}$ and $\widehat{\log S_i}^{(t)}$ through \eqref{em_s1} and \eqref{em_s3} using the updated parameters $\hat{\theta}_0^{(t)}$, $\hat{\sigma}_0^{(t)}$ and $\hat{\mu}^{(t)}$ as well as $\hat{\alpha}^{(t-1)}$.
    \item {\bf CM-step 2:} Compute $\hat{\alpha}^{(t)}$ with a numerical optimisation method using the updated expectations $\widehat{S_i}^{(t)}$ and $\widehat{\log S_i}^{(t)}$.
\end{itemize}
The algorithm encounters difficulties during the E-step when $X_i$ is very close to $\mu$. In this case, it is recommended to replace $\delta_i$ by $\delta_i^{*}=\Delta/\sqrt{2\alpha+\theta_0^2/\sigma_0^2}$ for some $\Delta>0$. Recommendations for values of $\Delta$ are given in Bee et al.\ \cite{bee18}. An improved version of the algorithm in the multivariate case
is available in Nitithumbundit and Chan \cite{niti19}.

\subsection{Further estimation techniques 
}
For the sake of completeness, we mention several other estimation methods in the literature. In Loregian, Mercuri and Rroji \cite{loreg12}, a Gaussian quadrature is used to approximate an integral representation of the VG density by a mixture of Gaussian densities, and then an ECM algorithm with respect to the mixing variable is used to estimate the parameters. Nzokem \cite{nzo21} uses a fractional Fourier transform  in order to approximate the VG density. In Finaly and Seneta \cite{fs08}, a minimum $\chi^2$-approach is studied, which compares theoretical and empirical quantiles (see also Finlay and Seneta \cite{fs06} and Tjetjep and Seneta \cite{ts06}). Finlay and Seneta \cite{fs08} further present an estimation procedure based on the characteristic function. It is worth noting that this approach is still valid when there is a certain dependency structure in the data. Finlay and Seneta \cite{fs08} describe a Bayesian procedure, and a Markov Chain Monte Carlo method is used in order to sample from the posterior distribution.  We also refer to Rathgeber, Johannes and  St\"{o}ckl \cite{rss16}, in which a number of the estimation techniques described in this section are compared on the basis of real financial data. 


\section{Applications}\label{sec4}

\subsection{Exact distribution in mathematical statistics 
}\label{sec4.1}




\subsubsection{Connection to sample correlation}

Let $(X_i,Y_i)$, $1\leq i\leq n$, be i.i.d.\ bivariate normal random vectors with mean vector $(\mu_X,\mu_Y)$, variances $(\sigma_X^2,\sigma_Y^2)$ and correlation coefficient $\rho$. It was shown by Pearson, Jefferey and Elderton \cite{p29} that the product-moment coefficient
\begin{equation*}p_{n}=\frac{1}{n}\sum_{i=1}^n(X_i-\overline{X})(Y_i-\overline{Y})
\end{equation*}
has a VG distribution (they referred to it as the Bessel function distribution). An alternative proof is given by Kotz et al.\ \cite[Proposition 4.1.5]{kkp01}. In our parametrisation, 
\begin{equation*}p_{n}\sim\mathrm{VG}(n-1,\rho\sigma_X\sigma_Y/n,\sigma_X\sigma_Y\sqrt{1-\rho^2}/n,0).
\end{equation*}
Here we applied the scaling property (\ref{muc})
 of the VG distribution 
to deal with the case of general variances $\sigma_X^2,\sigma_Y^2>0$; Kotz et al.\ \cite{kkp01} proved the result for
 the case 
$\sigma_X=\sigma_Y=1$. As observed by Kotz et al.\ \cite{kkp01}, when $n=3$ the product-moment coefficient $p_{3}$ has an asymmetric Laplace distribution (see part 1 of Section \ref{subcalss}). Moreover, as in part 4 of Section \ref{subcalss}, let $Z=UV$, where $(U,V)$ is a bivariate normal random vector with zero mean vector, variances $(\sigma_X^2,\sigma_Y^2)$ and correlation coefficient $\rho$. Consider also the mean $\overline{Z}_{n-1}=(n-1)^{-1}(Z_1+Z_2+\cdots+Z_{n-1})$, where $Z_1,Z_2,\ldots,Z_{n-1}$ are independent copies of $Z$. Then, by the relation (\ref{vgrep}) and the scaling property (\ref{muc}), we have that, for $n\geq2$,
\begin{align}\label{pn}p_{n}=_d (1-1/n)\overline{Z}_{n-1}.
\end{align}
In particular, $p_{2}=_dZ/2$. It follows from (\ref{pn}) and the central limit theorem that $\sqrt{n}(p_n-\rho\sigma_X\sigma_Y)\rightarrow_d N(0,\sigma_X^2\sigma_Y^2(1+\rho^2))$, as $n\rightarrow\infty$.

\subsubsection{Connection to the Wishart distribution}

Let $N$ be a $p\times n$ matrix, whose $n$ columns are independent $p$-dimensional multivariate normal random vectors with zero mean vector and positive definite covariance matrix $V=(v_{ij})$. Then the Wishart distribution
is the probability distribution of the $p\times p$ random matrix $X=NN^\intercal$. We denote such a random matrix by $W_p(V,n)$. The positive integer $n\geq1$ is referred to as the number of degrees of freedom, and when $V=p=1$, the Wishart distribution reduces to the $\chi_n^2$ distribution. The Wishart distribution is widely used in multivariate statistics, Bayesian analysis and random matrix theory.

A fundamental distributional property of the Wishart distribution is the marginal distribution of its entries. It is well-known that the diagonal entries are suitably scaled chi-square random variables: $X_{ii}\sim v_{ii} \chi_n^2$, $1\leq i\leq p$ (see, for example, Kollo \cite[Corollary 2.4.2.2]{kollo05}). It is perhaps less well-known that the off-diagonal entries are VG distributed (see Pearson et al.\ \cite{p29} for the two-dimensional case). More precisely, 
\begin{equation}
\label{xij}
X_{ij}\sim\mathrm{VG}(n,v_{ij},(v_{ii}v_{jj}-v_{ij}^2)^{1/2},0), \quad i\not=j.
\end{equation}
As $V$ is a covariance matrix, it is natural to set $v_{ii}=\sigma_{i}^2$ and $v_{ij}=\rho_{ij}\sigma_i\sigma_j$, $i\not=j$, where $\sigma_i^2$ is the variance of the $i$-th component of a $p$-dimensional $N_p(\mathbf{0},V)$ multivariate normal random vector, 
whilst $\rho_{ij}$ is the correlation coefficient of the $i$-th and $j$-th components.
 With this notation, we have that $X_{ij}\sim\mathrm{VG}(n,\rho_{ij}\sigma_i\sigma_j,\sigma_i\sigma_j(1-\rho_{ij}^2)^{1/2},0)$, $i\not=j$. From (\ref{vgrep}), it therefore follows that $X_{ij}$ has a neat representation as a sum of independent copies of correlated zero mean normal random variables. Let $Z_1,\ldots,Z_n$ be independent $\mathrm{PN}(\rho_{ij},\sigma_i,\sigma_j)$ random variables (recall that this notation was introduced in part 4 of Section \ref{subcalss}). Then
\begin{equation*}X_{ij}=_d Z_1+\cdots+Z_n, \quad i\not=j.
\end{equation*}

That the off diagonal entries are VG distributed can be deduced from the Bartlett decomposition of the Wishart distribution (see Anderson \cite{a03}) and the representation (\ref{rep2}) of the VG distribution in terms of independent normal and chi-square random variables. The Bartlett decomposition of $X\sim W_p(V,n)$ is 
\begin{equation}\label{bart}X=_d LAA^\intercal L^\intercal,
\end{equation}
where $L$ is the Cholesky factor of $V$ and $A$ is a lower triangular matrix with diagonal entries $S_i\sim \chi_{n-i+1}$ (so that $S_i^2\sim \chi_{n-i+1}^2$), $1\leq i\leq p$, and off-diagonal entries $0$ for $j<i$, and $N_{ij}\sim N(0,1)$ for $i>j$. All entries in the matrix $A$ are mutually independent. From the Bartlett decomposition (\ref{bart}) we obtain that
\begin{align}\label{x12}X_{12}=_d v_{12}S_1^2+(v_{11}v_{22}-v_{12}^2)^{1/2}S_1 N_{12},
\end{align}
where we used that $L_{11}=\sqrt{v_{11}}$, $L_{12}=v_{12}/\sqrt{v_{11}}$ and $L_{22}=(v_{22}-v_{12}^2/v_{11})^{1/2}$. It therefore follows from the representation (\ref{rep2}) of the VG distribution that $X_{12}\sim\mathrm{VG}(n,v_{12},(v_{11}v_{22}-v_{12}^2)^{1/2},0)$. We can use the Bartlett decomposition (\ref{bart}) to obtain similar representations of the off-diagonal elements $X_{ij}$, $i\not=j$, in terms of independent normal and chi-square random variables, although these representations are more complicated and it is harder to infer that the elements are indeed VG distributed through this approach. However, the fact that the off-diagonal entries are VG distributed with parameters given as in (\ref{xij}) follows from (\ref{x12}) and symmetry considerations. This can be seen analytically, for a general off-diagonal element $X_{ij}$, by suitably reordering the elements in the Wishart matrix $X$ and making the corresponding reordering in the matrix $V$, and then obtaining the Bartlett decomposition of the new Wishart matrix $X'$. The entry $X_{ij}$ (entry $X_{12}'$ in the re-ordered Wishart matrix $X'$) then can be seen to have a representation of the form $X_{ij}=_d v_{ij}S_1^2+(v_{ii}v_{jj}-v_{ij}^2)^{1/2}S_1 N_{12}$, and the claim follows.


\subsection{Variance-gamma process}\label{sec4.2}

In Section \ref{sec2.3}, we saw that the VG distribution is infinitely divisible. As a consequence of the theory of infinitely divisible distributions and processes (see Ferguson and Klass \cite{fk72}), we can define a L\'evy process with VG distributed increments. The resulting L\'evy process is known as the VG process, which, as we shall discuss in Section \ref{sec4.3}, is widely used in financial modelling. The terminolgy \emph{Laplace motion} is also often employed; see, for example, Kotz et al.\ \cite{kkp01} and Kozubowski and Podg\'{o}rski \cite{kp12}. The first complete presentation of the VG model was given by Madan and Seneta \cite{madan} for the case of symmetric VG distributed increments ($\theta=0$), which was extended to incorporate skewness $(\theta\in\mathbb{R}$) by Madan and Milne \cite{mm91} and Madan et al.\ \cite{mcc98}. A few years before the seminal paper of Madan and Seneta \cite{madan}, the same authors published two papers (Madan and Seneta \cite{ms87a,ms87b}) which contained some basic properties of the VG process; see Seneta \cite{s07} for an account of the early years of the VG process.

In this section, we shall review some of the most fundamental properties of the VG process and briefly discuss why these properties are desirable in the context of financial modelling. 
We shall only consider the univariate VG process; for an account of the multivariate VG process, which was introduced by Madan and Seneta \cite{madan}, we refer the reader to Luciano and Schoutens \cite{ls06}, Semeraro \cite{sem08}, Luciano and Semeraro \cite{ls10} and Luciano, Marena and Semeraro \cite{lms16}, and Linders and Stassen \cite{ls15} for applications to basket options calibration on the DJIA index. We also do not discuss methods for simulating VG processes, and refer the reader to Avramidis and L'Ecyuer \cite{ae06}, Fu \cite{f07} and Korn, Korn and Kroisandt \cite[Section 7.3.3]{3k} for algorithms.

Madan et al.\ \cite{mcc98} constructed the VG process as a subordinated Brownian motion with a time change by a gamma process. Consider a Brownian motion with drift $\theta$ and volatility $\sigma$ given by $b(t;\theta,\sigma)=\theta t+\sigma B(t)$, where $B(t)$ is standard Brownian motion. Recall that the gamma process $\gamma(t;\mu,\nu)$ with mean rate $\mu$ and variance rate $\nu$ is a L\'evy process with increment $\tau_h=\gamma(t+h;\mu,\nu)-\gamma(t;\mu,\nu)$ following the $\Gamma(\mu^2h/\nu,\mu/\nu)$ distribution. Suppose that the processes $B(t)$ and $\gamma(t;\mu,\nu)$ are independent. Then, the VG process $X(t)=X(t;\sigma,\nu,\theta)$, with parameters $\nu>0$, $\theta\in\mathbb{R}$ and $\sigma>0$, is defined as
\begin{equation}\label{vgp}X(t)=b(\gamma(t;1,\nu);\theta,\sigma).
\end{equation}
As we shall see, control over skewness and kurtosis can be attained via the parameters $\theta$ and $\nu$, respectively. This feature
is important in financial modelling. As noted by Madan and Seneta \cite[p.\ 517]{madan}, the random gamma time change in (\ref{vgp}) has an interesting economic interpretation of supposing that economically relevant time is random in that the market (price process $X(t)$) moves at different speeds on different days, and the process $\gamma(t;1,\nu)$ is a measure of this speed. Samples paths showing the effect of varying the parameter $\nu$ in the symmetric VG process are given in Figure \ref{fig:vgprocess}.

\begin{figure}[!]
    \centering
\vspace{.2cm}
   \begin{subfigure}{.49\textwidth}
\captionsetup{width=.95\textwidth}
  \centering
  \includegraphics[width=7.4cm]{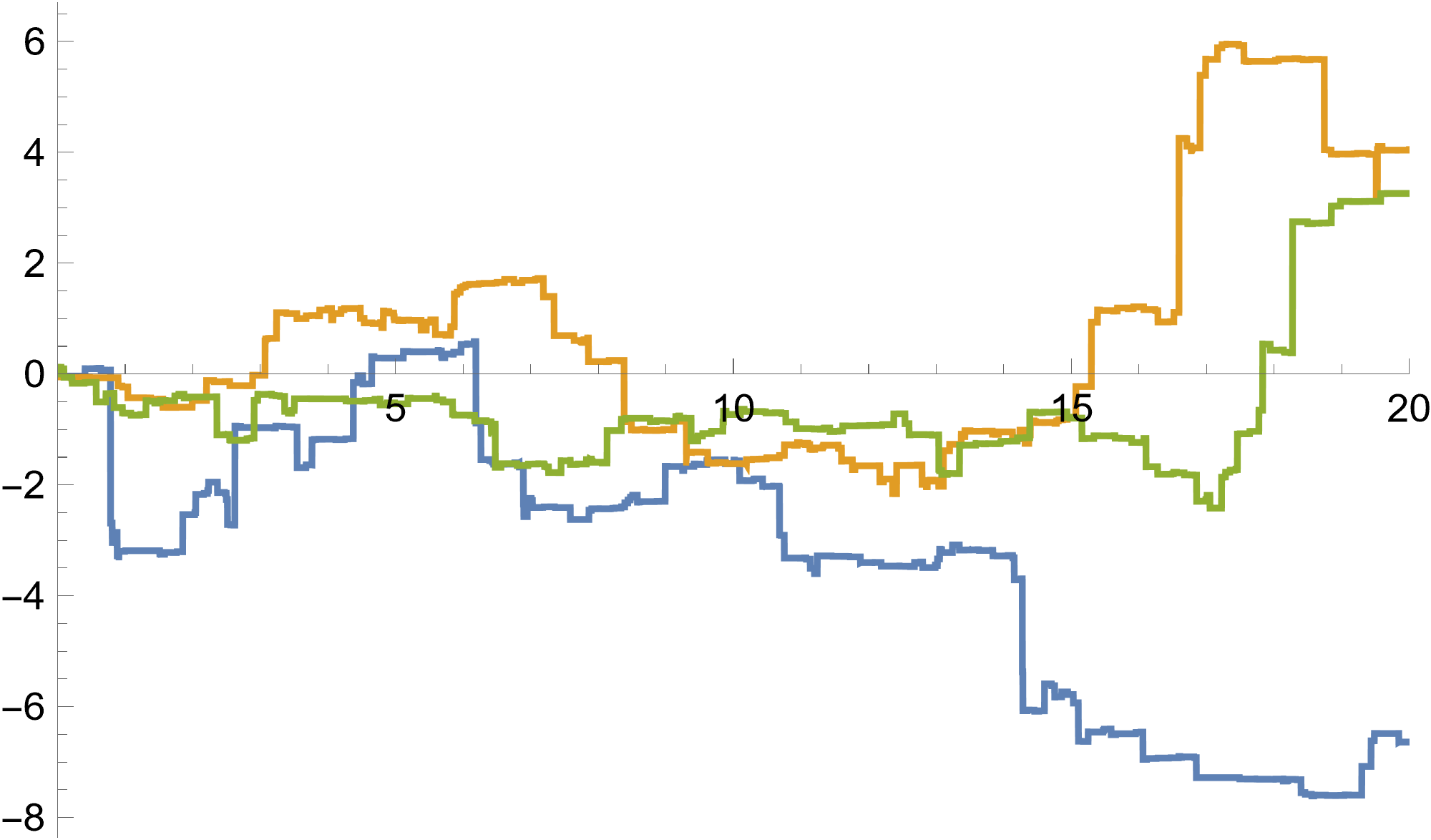}
\end{subfigure}%
  \begin{subfigure}{.49\textwidth}
\captionsetup{width=.95\textwidth}
  \centering
  \includegraphics[width=7.4cm]{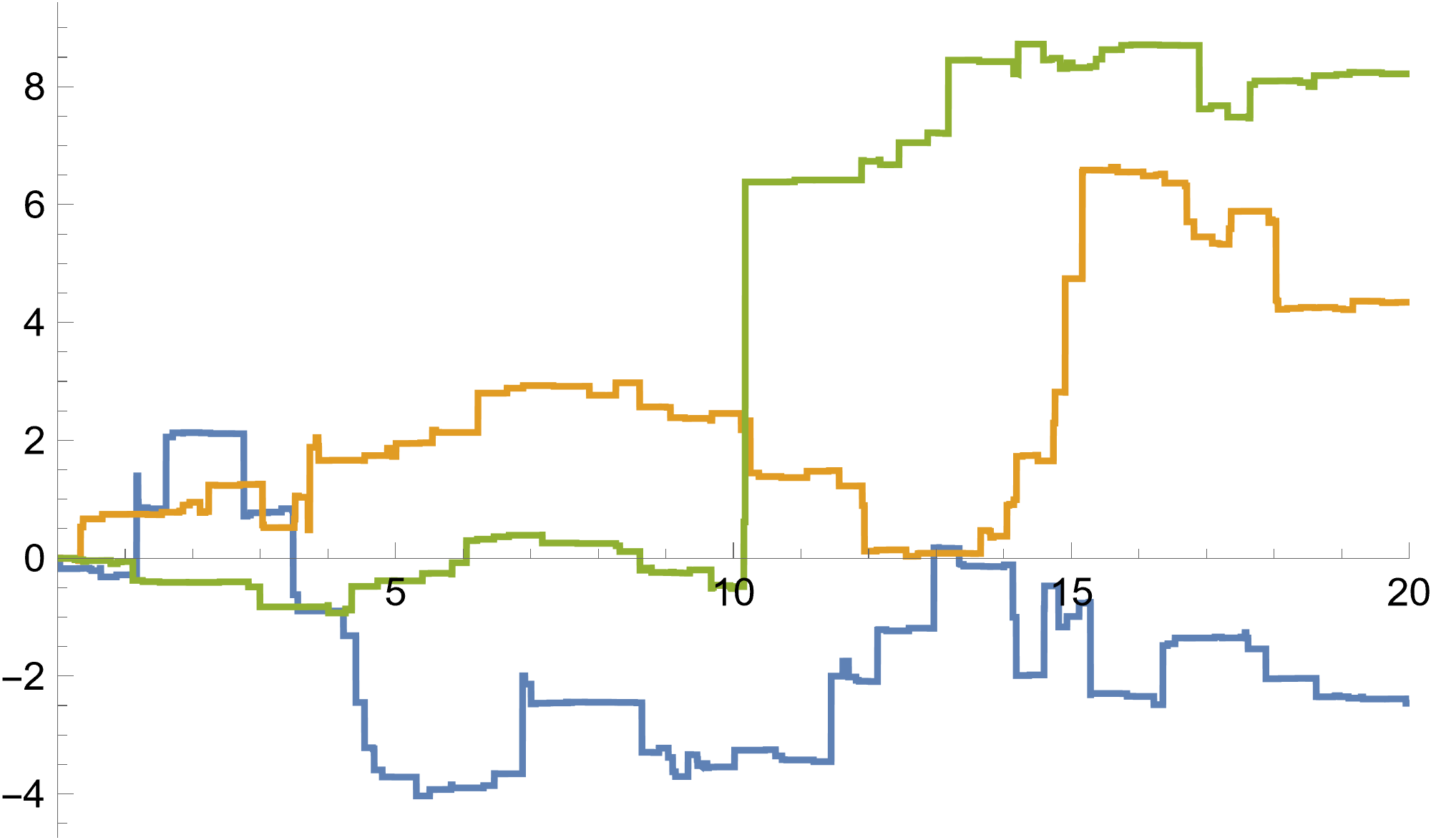}
\end{subfigure}
\caption{\label{fig:vgprocess}\it Three realisations of the VG process $X(t;1,1,0)$ (left image) and  $X(t;1,2,0)$ (right image).}
\end{figure}

To see that the increments of the VG process (\ref{vgp}) are VG distributed, we observe that 
\begin{align*}X(t+h)-X(t)&= \theta\tau_h+\sigma(B(\gamma(t+h;1,\nu))-B(\gamma(t;1,\nu))=_d\theta\tau_h+\sigma\tau_h^{1/2}B(1).
\end{align*}
As $B(1)\sim N(0,1)$ and $\tau_h\sim \Gamma(\mu^2h/\nu,\mu/\nu)$ are independent, it follows from 
(\ref{rep2}) that 
\begin{equation}\label{vgproch}X(t+h)-X(t)\sim\mathrm{VG}(2h/\nu,\theta\nu/2,\sigma\sqrt{\nu/2},0).
\end{equation}
In particular, since $X(0)=0$,
\begin{equation}\label{vgproc}X(t)\sim\mathrm{VG}(2t/\nu,\theta\nu/2,\sigma\sqrt{\nu/2},0).
\end{equation}


As we shall see shortly, the VG process is of finite variation, so can be expressed as the difference of two independent increasing processes. These increasing processes are in fact gamma processes, and we have the representation (see Madan et al.\ \cite{mcc98})
\begin{equation}\label{gamx}X(t)=_d\gamma_p(t;\mu_p,\nu_p)-\gamma_n(t;\mu_n,\nu_n),
\end{equation}
where
\begin{align*}\mu_p=\frac{1}{2}\sqrt{\theta^2+\frac{2\sigma^2}{\nu}}+\frac{\theta}{2}, \quad \mu_n=\frac{1}{2}\sqrt{\theta^2+\frac{2\sigma^2}{\nu}}-\frac{\theta}{2}, \quad \nu_p=\mu_p^2\nu, \quad \nu_n=\mu_n^2\nu.
\end{align*}
Here $\gamma_p(t;\mu_p,\nu_p)$ plays the role of a returns process, whilst $\gamma_n(t;\mu_n,\nu_n)$ is a losses process. The representation (\ref{gamx}) of the VG process as a difference of two independent gamma processes follows from the representation (\ref{gamrep}) and the fact that $\gamma(t;\mu,\nu)\sim \Gamma(\mu^2t/\nu,\mu/\nu)$. 

From (\ref{vgproc}) and the moments formulas (\ref{mean}) and (\ref{mom2})--(\ref{mom4}) we obtain that
\begin{align*}\mathbb{E}[X(t)]&=\theta t, \\
\mathbb{E}[(X(t)-\mathbb{E}[X(t)])^2]&=(\theta^2\nu+\sigma^2)t, \\
\mathbb{E}[(X(t)-\mathbb{E}[X(t)])^3]&=(2\theta^3\nu^2+3\sigma^2\theta\nu)t, \\
\mathbb{E}[(X(t)-\mathbb{E}[X(t)])^4]&=(3\sigma^4\nu+12\sigma^2\theta^2\nu^2+6\theta^4\nu^3)t+(3\sigma^4+6\sigma^2\theta^2\nu+3\theta^4\nu^2)t^2
\end{align*}
(see Madan et al.\ \cite{mcc98}). We see that the parameter $\theta$ does not directly control skewness, although skewness is zero when $\theta=0$, and positive when $\theta>0$.
 When $\theta=0$, the kurtosis is 
$3(1+\nu/t)$, and so $\nu$ is the percentage of excess kurtosis for a unit interval of time $t=1$. As $t$ increases, the kurtosis decreases to 3, the kurtosis of the standard normal distribution. As noted by Madan and Seneta \cite{madan}, this is consistent with the empirical evidence that daily returns are heavy tailed, whilst monthly returns are normally distributed.

From (\ref{vgproc}) and (\ref{cfcf}), we have that the characteristic function of the VG process is given by
\begin{equation}\label{cfcfproc}\varphi_{X(t)}(u)=\mathbb{E}[\mathrm{e}^{\mathrm{i}uX(t)}]=\big(1-\mathrm{i}\theta \nu u+(\sigma^2\nu/2)u^2\big)^{-t/\nu}
\end{equation}
(see Madan et al.\ \cite{mcc98}). That the characteristic function (\ref{cfcfproc}) takes a simple form is important in option pricing, since, for example, the fast fourier transform approach from the classic paper of Carr and Madan \cite{cm99} can be applied to numerically determine option values. From (\ref{levyrep}) and (\ref{vgproc}) we have that the L\'evy-Hinchin representation of the VG process is given by $\varphi_{X(t)}(u)=\exp\big(t\int_\mathbb{R}(\mathrm{e}^{\mathrm{i}ux}-1)\nu_{VG}(x)\,\mathrm{d}x\big)$, with L\'evy density
\begin{equation*}\nu_{VG}(x)=\frac{\mu_n^2}{\nu_n|x|}\mathrm{e}^{-\mu_n|x|/\nu_n}\mathbf{1}_{x<0}+\frac{\mu_p^2}{\nu_p x}\mathrm{e}^{-\mu_px/\nu_p}\mathbf{1}_{x>0}.
\end{equation*}
From its L\'evy-Hinchin representation, it follows that the VG process is a pure jump process (there is no diffusion component). 
The pure jump nature of the VG process represented a departure from much of the literature on option pricing; for example, the Black-Scholes model is a pure diffusion model. It was noted by Bakshi, Chen and Cao \cite{bcc} that a jump component is important in option pricing, as pure diffusion models have difficulties in explaining smile effects, in particular in short-dated option prices. Moreover, as asserted by Madan et al.\ \cite{mcc98}, the Black-Scholes model is a parametric special case of the VG model (recall from Section \ref{subcalss} that the normal distribution is a limiting case of the VG distribution), meaning high activity is already accounted for and so it is not necessary to include a diffusion component in addition.
Indeed, due to the factor of $|x|^{-1}$, the L\'evy density integrates to infinity, and so the VG process is of infinite activity. Also, since $|x|$ is integrable with respect to the L\'evy density, the process is of finite variation. As noted by Carr et al.\ \cite{cgmy}, finite variation processes may be more useful than infinite variation processes in explaining the measure change from the statistical to the risk neutral process, since they allow more flexibility between the local characteristic of the martingale components under the two measures. Finally, again as noted by Carr et al.\ \cite{cgmy}, the VG process has a completely monotone L\'evy density. In particular, the derivative of the L\'evy density is positive for negative jumps and negative for positive jumps, so that large jumps arrive less frequently than small jumps, which is consistent with what we observe in price movements. Indeed, most jumps are of a vanishingly small size; for a compound Poisson approximation and a detailed account of the fine structure of the jumps see Madan and Seneta \cite{madan} and Kotz et al.\ \cite[Section 4.2]{kkp01}.



\subsection{Financial modelling}\label{sec4.3}

The Black-Scholes model is the paradigm model in mathematical finance. Under this model, the asset price $S(t)$ at time $t$ obeys geometric Brownian motion
\begin{equation}\label{gbm}S(t)=S(0)\exp((\mu-\sigma^2/2) t+\sigma B(t)),
\end{equation}
where $\mu\in\mathbb{R}$ and $\sigma>0$ are constants denoting drift and volatility, and $B(t)$ is standard Brownian motion. Whilst widely used in financial modelling, a number of deficiencies have been noted. For example, as detailed by Heyde \cite{h99}, under geometric Brownian motion, the log returns $Y(t)=\log S(t)-\log S(t-1)$ are i.i.d.\ normal random variables, which is often in contrast to empirical log returns data; see Eberlein and Keller \cite{eberlein1} and Leonenko, Petherick and Sikorskii \cite{leo11} for empirical studies. In particular, the tails of the normal distribution are not heavy enough for modelling of daily returns; this was one of the original motivations of Madan and Seneta \cite{madan} to introduce the VG process with an additional parameter to control kurtosis to allow for heavier tails.
On account of the flexibility to control skewness and kurtosis to allow for a better empirical fit to financial data, together with other desirable properties such as those discussed in Section \ref{sec4.2}, the VG model has become widely used in financial modelling; a comparison of the fit of the VG distribution and the normal distribution to real financial data is given in Figure \ref{fig:finance_application}. In this section, we give a brief introduction to
an extensive literature.

\begin{figure}[!]
    \centering
\vspace{.2cm}
   \begin{subfigure}{.49\textwidth}
\captionsetup{width=.95\textwidth}
  \centering
  \includegraphics[width=7.4cm]{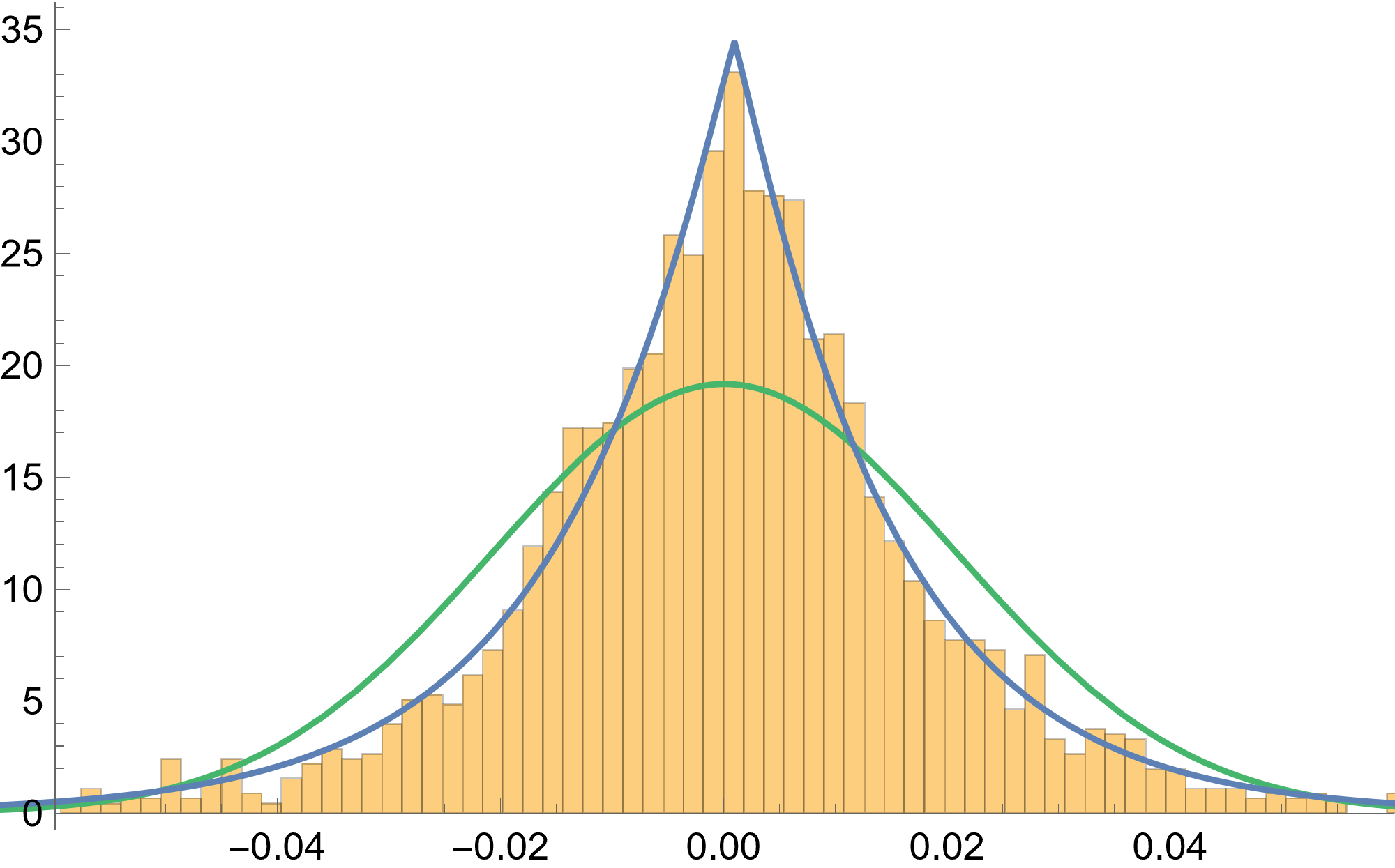}
\end{subfigure}%
  \begin{subfigure}{.49\textwidth}
\captionsetup{width=.95\textwidth}
  \centering
  \includegraphics[width=7.4cm]{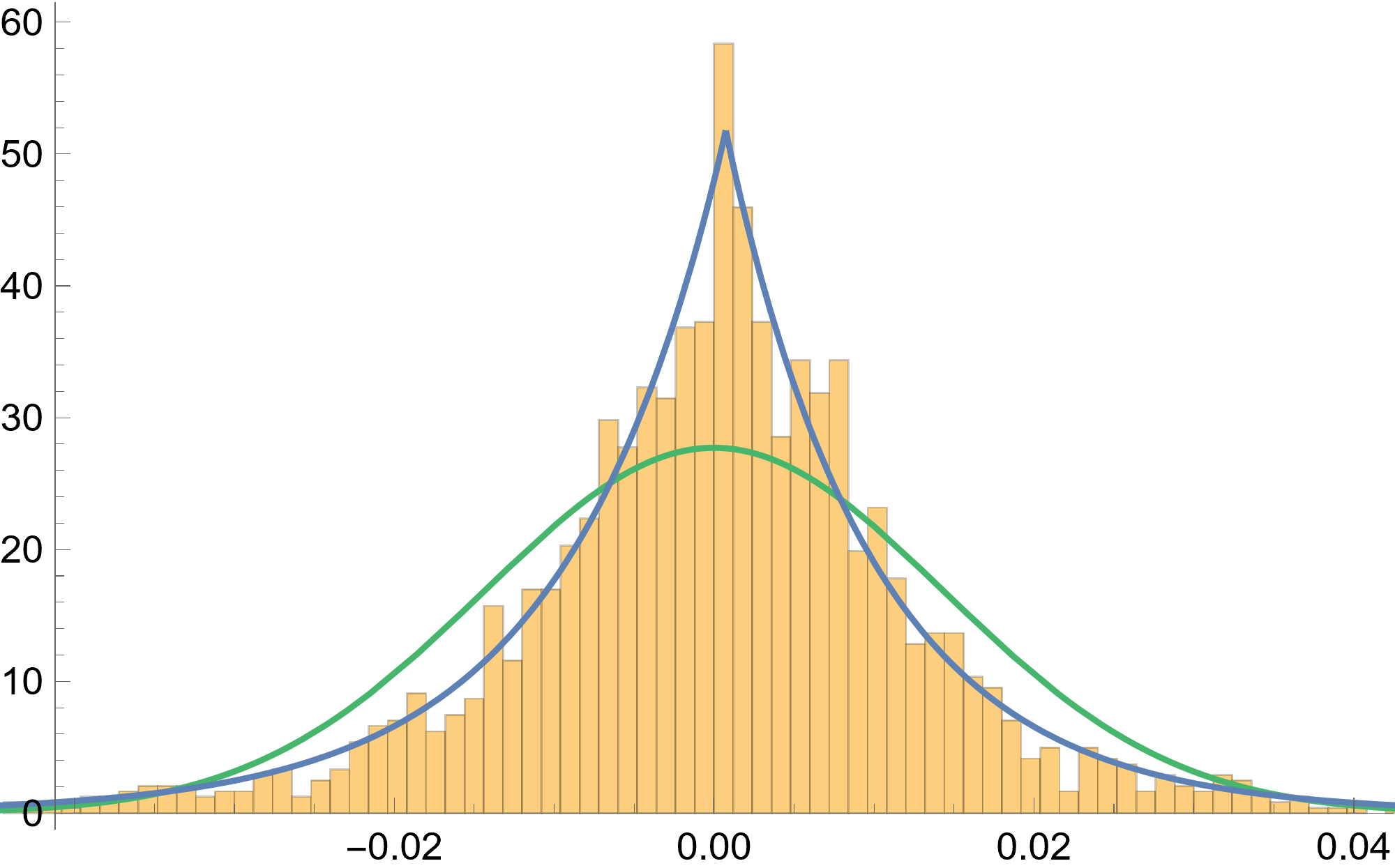}
\end{subfigure}
\caption{\label{fig:finance_application}\it A real data example: We see the daily log returns of the Mercedes Benz Group AG (left picture, 01.01.2013 -- 01.01.2023, 2517 observations) and the Global X DAX Germany ETF (right picture, 01.01.2015 -- 01.01.2023, 2013 observations) fitted to a VG distribution (\includegraphics[scale=1]{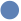}) and a normal distribution (\includegraphics[scale=1]{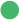}). We used maximum likelihood estimation in order to fit the VG distribution.
} 
\end{figure}

The statistical and risk neutral dynamics of a stock price in terms of the VG process was given by Madan et al.\ \cite{mcc98}. The statistical stock price is obtained by replacing the role of Brownian motion in the geometric Brownian motion model (\ref{gbm}) by the VG process:
\begin{equation}\label{expvg}S(t)=S(0)\exp(mt+X(t;\sigma_S,\nu_S,\theta_S)+\omega_St),
\end{equation}
where $X(\cdot)$ is the VG process, $m$ is the mean rate of return on the stock under the statistical probability measure and the subscript $S$ in the VG parameters stresses that they are the statistical parameters. The value $\omega_S=\nu_S^{-1}\log(1-\theta_S\nu_S-\sigma_S^2\nu_S/2)$ is determined by the non-arbitrage condition to ensure that $\mathbb{E}[S(t)]=S(0)\exp(mt)$, that is $\mathbb{E}[\exp(X(t))]=\exp(-\omega_St)$, with $\mathbb{E}[\exp(X(t))]$ evaluated using (\ref{vgproc}) and the moment generating function formula (\ref{mgf}).

Under the risk neutral process, stock prices discounted at the risk free interest rate are martingales. Therefore the mean rate of return on the stock under this probability measure is the continuously compounded interest rate $r$. The risk neutral process is
\begin{equation}\label{rnvg}S(t)=S(0)\exp(rt+X(t;\sigma_{RN},\nu_{RN},\theta_{RN})+\omega_{RN}t),
\end{equation}
where the subscripts $RN$ indicate that these are the risk free parameters, and $\omega_{RN}=\nu_{RN}^{-1}\log(1-\theta_{RN}\nu_{RN}-\sigma_{RN}^2\nu_{RN}/2)$, derived using the same considerations as for the statistical stock price. As discussed in detail by Madan et al.\ \cite{madan}, unlike for diffusion based price processes, the statistical and risk free parameters do not need to be equal, and are often quite different.

As noted by Madan et al.\ \cite{madan}, the log stock price relative $\log(S(t)/S(0))$ when prices follow the VG process dynamics of (\ref{expvg}) is VG distributed (this follows from (\ref{muc}) and (\ref{vgproc})):
\begin{equation*}\log(S(t)/S(0))\sim \mathrm{VG}(2t/\nu,\theta\nu/2,\sigma\sqrt{\nu/2},(m+\omega)t),
\end{equation*}
where $\omega=\nu^{-1}\log(1-\theta\nu-\sigma^2\nu/2)$ and we have dropped the subscript $S$ from the parameters. This time using (\ref{vgproch}), we see that the log returns are also VG distributed: 
\begin{equation*}Y(t)=\log S(t)-\log S(t-1)\sim \mathrm{VG}(2/\nu,\theta\nu/2,\sigma\sqrt{\nu/2},m+\omega);
\end{equation*}
see Finlay and Seneta \cite{fs06}, who also give formulas to describe the covariance structure of the log returns process $\{Y(t)\}$.

Madan et al.\ \cite{madan} derived a closed-form formula for the price of a European call option when the risk neutral dynamics of the stock price follows the VG model (\ref{rnvg}). A standard result states that the price of a European call option $c(S(0);K,t)$, for a strike $K$ and maturity $t$, is given by
\begin{equation}\label{exmax}c(S(0);K,t)=\mathrm{e}^{-rt}\mathbb{E}[\mathrm{max}(S(t)-K,0)].
\end{equation}
Building on the approach of Madan and Milne \cite{mm91}, Madan et al.\ \cite{mcc98} evaluated the expectation in (\ref{exmax}) by first conditioning on the random gamma process time change, under which the VG process $X(t)$ is normally distributed (see part 1 of Section \ref{sec2.4}) and the option value is given by a Black-Scholes type formula. Finally, the European option price for VG risk neutral dynamics is attained on integrating the resulting conditional Black-Scholes formula with respect to the gamma density. The formula is
\begin{align}c(S(0);K,t)&=S(0)\Psi\bigg(d\sqrt{\frac{1-c_1}{\nu}},(\alpha+s)\sqrt{\frac{\nu}{1-c_1}},\frac{t}{\nu}\bigg)\nonumber\\
\label{vgbs}&\quad-K\mathrm{e}^{-rt}\Psi\bigg(d\sqrt{\frac{1-c_2}{\nu}},\alpha s\sqrt{\frac{\nu}{1-c_2}},\frac{t}{\nu}\bigg),
\end{align}
where $\alpha=-\theta/\sqrt{\sigma^2+\theta^2\nu/2}$, $c_1=\nu(\alpha+s)^2/2$, $c_2=\nu\alpha^2/2$, and
\begin{align*}d=\frac{1}{s}\bigg[\log\bigg(\frac{S(0)}{K}\bigg)+rt+\frac{t}{\nu}\log\bigg(\frac{1-c_1}{1-c_2}\bigg)\bigg],
\end{align*}
and $\Psi$ is expressed explicitly in terms of the modified Bessel function of the second kind and the degenerate hypergeometric function of two variables; see
Madan et al.\ \cite{mcc98} for the formula. As noted by Madan et al.\ \cite{mcc98}, on letting $\nu\rightarrow0$ in (\ref{vgbs}), we recover the classical Black-Scholes option pricing formula. 

Madan et al.\ \cite{mcc98} compared the pricing performance on S\&P 500 option data of the pricing formula (\ref{vgbs}) with the Black-Scholes model and the symmetric special case ($\theta=0$). They found that the performance was superior over the Black-Scholes model, and provided evidence that whilst the statistical density of the underlying risk is typically symmetric, the risk neutral density implied option data is negatively skewed, demonstrating the importance of considering the general non-symmetric VG model.

Beyond the work of Madan et al.\ \cite{mcc98}, 
the VG model has been successfully tested on real market data and has been found to outperform that Black-Scholes model and Jump-Diffusion models in a number of other situations, such as European-style options on the HSI index by Lam, Chang and Lee \cite{lcl02} and currency options by Daal and Madan \cite{dm05}. Closed-form pricing formulas have also been obtained for various path-independent payoffs under the VG model; see, for example, Aguilar \cite{a20} and Ivanov \cite{i18}. 
Typically, however, closed-form option formulas are unavailable; see, for example, Hirsa and Madan \cite{hm04}, in which a partial integro-differential equation is derived for pricing American options under the VG process. In such cases, numerical methods are required, and we refer the reader to the thesis of Fiorani \cite{f04} for an extensive account of methods, as well as the recent article of Aguilar \cite{a20} for an overview of some popular numerical methods.



\subsection{Time series modelling}\label{sec4.4}

Recall from Section \ref{sec4.3}, that the VG distribution is a useful alternative to the normal distribution in financial modelling on account of its heavier tails and capacity to incorporate skewness. This applies equally well to innovations in time series modelling. As discussed in detail by Johannesson \cite{ar1}, VG distributed innovations are a natural alternative to the widely used normally distributed innovations, on account of a natural interpretation of conditionally normally distributed innovations with gamma distributed variance (which results in heavier tails), whilst retaining some analytic tractability and are relatively simple to simulate from. In this section, we provide a brief overview of a growing literature on time series models with VG distributions for innovations. 

 In Tomy and Jose \cite{gnl-ar1}, an autoregression of order $1$ with generalized Gaussian-Laplace innovations was studied:
\begin{equation}
\label{eq:ar1}
X_t = aX_{t-1} + Z_t + V_t,
\end{equation}
where $a$ is a model parameter and $Z_t \sim  N(0, \rho^2)$ and $ V_t \sim \mathrm{VG}(r, \theta, \sigma, \mu)$ are independent.
Here the distribution for the innovations is more general than the VG distribution: it is the sum of independent VG and normal random variables (centered, without loss of generality). Tomy and Jose \cite{gnl-ar1} fitted parameters using the classic method of moments, by computing the first six moments; however, the problem of goodness-of-fit with this method was not discussed. Nitithumbundit and Chan \cite{niti19} considered autoregression of order $1$ in a slightly different setting than in \eqref{eq:ar1}. They only had a VG component $V_t$ (no $Z_t$), but the setting was multivariate. The authors successfully used ECM (expectation-conditional maximisation) to fit the model even in the case of unbounded VG density of innovations. They applied their methods to financial data, and reported a superior fit to the data than analogous models with multivariate normally distributed innovations. More general ARMA($p,q$) models driven by asymmetric Laplace (the $\mathrm{VG}(2,\theta,\sigma,\mu)$ distribution) noise and ARMA($p,q$) models driven by GARCH asymmetric Laplace noise were studied by Trindade, Zhu and Andrews \cite{time-series}. Similarly to Nitithumbundit and Chan \cite{niti19}, they applied their models to financial data, again reporting a better fit to the data than for corresponding models with normal innovations.




Johannesson \cite{ar1} considered a more complicated time series model in the context of road topography. For innovations, a symmetric VG random variable $Y =\mu+ \sigma\sqrt{S}T $ is taken from \eqref{rep2} (with asymmetry parameter $\theta = 0$), replacing $T$ with 
$T_t$, an autoregressive process of order 1 with Gaussian innovations, and $S$ with $S_t$, a (nonlinear) autoregressive process of order 1 constructed so that $S_t$ has a gamma distribution for each $t$ and is stationary. They presented two methods to fit the autoregressive parameters, and fitted the model to real-life road topology data, showing that it satisfactorily retrieves the distributions of the length, average slope and height of hills. 

Finally, we mention that there has been recent interest in the study of non-Gaussian random fields built on the VG distribution, and we refer the reader to {\AA}berg, Podg\'orski and Rychlik \cite{apr09}, {\AA}berg and Podg\'orski \cite{ap11}, Bolin \cite{bolin}, Baxevani, Podg\'orski and Wegener \cite{bpw14} and the survey of Kozubowski and Podg\'{o}rski \cite{kp12} for an overview, as well as Podg\'{o}rski and Wegener \cite{pw11} for a discussion of estimation methods, and Bogsj\"o,  Podg\'orski and Rychlik \cite{bpr} for an application to the modelling of road surface irregularities.


\subsection{Approximation on Wiener space}\label{sec4.5}

In Section \ref{sec2.3}, we saw that the VG distribution is infinitely divisible and in Section \ref{sec2.4} we saw that the VG distribution has simple representations in terms of independent standard normal and gamma random variables. In virtue of these properties, the VG distribution is a natural candidate as a limiting distribution. 
Indeed, there has been recent interest in VG approximation on Wiener space. VG approximations for double Wiener-It\^{o} integrals have been studied by Azmoodeh, Eichelsbacher and Th\"ale \cite{aet21}, Eichelsbacher and Th\"ale \cite{eichelsbacher}
 and Gaunt \cite{gaunt vg3}
via Stein's method for VG approximation, which was first developed by Gaunt \cite{gaunt vg}, for which the starting point is the Stein characterisation (\ref{char1}).

\subsubsection{A six moment theorem for double Wiener-It\^{o} integrals}

Let $\mathfrak{H}$ be a real separable Hilbert space and denote by $\mathfrak{H}^{\odot 2}$ the second symmetric tensor product of $\mathfrak{H}$. For $f\in \mathfrak{H}^{\odot 2}$, we denote the double Wiener-It\^{o} integral by $I_2(f)$ (see Nourdin and Peccati \cite[Section 2.7]{np12} for a definition and fundamental properties). If $f\in L^2([0,T]^2,\mathrm{d}t)$ is symmetric then
\[I_2(f)=\int_{[0,T]^2}f(t_1,t_2)\,\mathrm{d}B_{t_1}\,\mathrm{d}B_{t_2},\]
where $B=(B_t)_{t\in[0,T]}$ is a standard two-dimensional Brownian motion (see Nourdin and Peccati \cite[Exercise 2.7.6]{np12}). Consider 
the Wasserstein $d_{\mathrm{W}}(F,G)$ and smooth Wasserstein distance $d_{2}(F,G)$ between the distributions of two random elements $F$ and $G$, defined by
\begin{align*}d_{\mathcal{H}}(F,G):&=\sup_{h\in\mathcal{H}}|\mathbb{E}[h(F)]-\mathbb{E}[h(G)]|,
\end{align*}
where the supremum is taken over the function classes
\begin{align*}
\mathcal{H}_{\mathrm{W}}&=\{h:\mathbb{R}\rightarrow\mathbb{R}\,|\,\text{$h'$ is Lipschitz, $\|h'\|_\infty\leq1$}\},\\\mathcal{H}_2&=\{h:\mathbb{R}\rightarrow\mathbb{R}\,|\,\text{$h'$ is Lipschitz, $\|h'\|_\infty\leq1$, $\|h''\|_\infty\leq1$}\},
\end{align*}
respectively, where the supremum norm of $g:\mathbb{R}\rightarrow\mathbb{R}$ is defined by $\|g\|_\infty=\sup_{x\in\mathbb{R}}|g(x)|$. Clearly, $d_{2}(F,G)\leq d_{\mathrm{W}}(F,G)$ for any random elements $F$ and $G$ for which $d_{\mathrm{W}}(F,G)$ is well-defined. 

Let $F_n=I_2(f_n)$ with $f_n\in \mathfrak{H}^{\odot 2}$, $n\geq1$. Following Eichelsbacher and Th\"ale \cite{eichelsbacher}, we write $\mathrm{VG}_c(r,\theta,\sigma)$ for $\mathrm{VG}(r,\theta,\sigma,-r\theta)$, the VG distribution with zero mean.  Then, Theorem 5.8 of Eichelsbacher and Th\"ale \cite{eichelsbacher} gives that, as $n\rightarrow\infty$, the sequence $(F_n)_{n\geq1}$ converges in distribution to $Y\sim\mathrm{VG}_c(r,\theta,\sigma)$ if and only if $\kappa_i(F_n)\rightarrow\kappa_i(Y)$, $i=2,3,4,5,6$. The cumulants $\kappa_i(Y)$ are readily obtained through the formulas of Section \ref{sec2.7} and the formula $\kappa_j(S+c)=\kappa_j(S)$, $j\geq2$, for any $c\in\mathbb{R}$. This ``six moment" theorem implies that the convergence of a sequence of double Wiener-It\^{o} integrals to the $\mathrm{VG}_c(r,\theta,\sigma)$ distribution is governed solely by the behaviour of the first six cumulants (equivalently, first six moments). This is a VG analogue of the celebrated ``fourth moment" theorem of Nualart and Peccati \cite{np05} for normal approximation of multiple Wiener-It\^{o} integrals.

Furthermore, quantitative ``six moment" theorems are also available. Let
\[\mathbf{M}(F_n)=\max\{|\kappa_i(F_n)-\kappa_i(Y)|\,:\,i=2,3,4,5,6\}.\]
Then, there exists a universal constant $C>0$, which only depends on $r$, $\theta$, $\sigma$, such that
\begin{equation}\label{optvg0}d_{\mathrm{W}}(F_n,Y)\leq C\sqrt{\mathbf{M}(F_n)}
\end{equation}
(Eichelsbacher and Th\"ale \cite{eichelsbacher} and Gaunt \cite{gaunt vg3}), and, in the weaker $d_2$ metric, it was shown by Azmoodeh et al.\ \cite{aet21} that  there exist universal constants $C_1,C_2>0$ such that
\begin{equation}\label{optvg}C_1\mathbf{M}(F_n)\leq d_{2}(F_n,Y)\leq C_2\mathbf{M}(F_n).
\end{equation} 
The rate of convergence in (\ref{optvg}) is optimal, and thus constitutes an analogue of the optimal fourth moment theorem of Nourdin and Peccati \cite{np15} for normal approximation.

\subsubsection{The generalized Rosenblatt process at extreme critical exponent}

We now review how Azmoodeh et al.\ \cite{aet21}
 and Gaunt \cite{gaunt vg3} used the bounds (\ref{optvg0}) and (\ref{optvg}) to derive bounds on the rate of convergence for a remarkable limit theorem of Bai and Taqqu \cite{bt17}. 
Consider the generalized Rosenblatt process $Z_{\gamma_1,\gamma_2}(t)$, which was introduced by Maejima and Tudor \cite{mt12} as the double Wiener-It\^{o} integral 
\begin{equation*}Z_{\gamma_1,\gamma_2}(t)=\int_{\mathbb{R}^2}^{\prime}\bigg(\int_0^t(s-x_1)_+^{\gamma_1}(s-x_2)_+^{\gamma_2}\,\mathrm{d}s\bigg)\,\mathrm{d}B_{x_1}\,\mathrm{d}B_{x_2},
\end{equation*}
where $B_{x}$ is standard Brownian motion, $\gamma_i\in(-1,-1/2)$, $i=1,2$, and $\gamma_1+\gamma_2>-3/2$, and the prime $\prime$ indicates exclusion of the diagonals $x_1=x_2$ in the stochastic integral. 
The Rosenblatt process is $Z_\gamma(t)=Z_{\gamma,\gamma}(t)$, $-3/4<\gamma<-1/2$ (see Taqqu \cite{t75}). As $Z_{\gamma_1,\gamma_2}(t)=_d t^{2+\gamma_1+\gamma_2}Z_{\gamma_1,\gamma_2}(1)$,
 we will consider the random variable $Z_{\gamma_1,\gamma_2}(1)$; results for general $t>0$ follow from rescaling.
For $\rho\in(0,1)$, define the random variable $Y_\rho$ by
\begin{equation*}Y_\rho=\frac{a_\rho}{\sqrt{2}}(X_1-1)-\frac{b_\rho}{\sqrt{2}}(X_2-1),
\end{equation*}
where $X_1$ and $X_2$ are independent $\chi_{1}^2$ random variables and 
\begin{align*}a_\rho=\frac{(2\sqrt{\rho})^{-1}+(\rho+1)^{-1}}{\sqrt{(2\rho)^{-1}+2(\rho+1)^{-2}}}, \quad b_\rho=\frac{(2\sqrt{\rho})^{-1}-(\rho+1)^{-1}}{\sqrt{(2\rho)^{-1}+2(\rho+1)^{-2}}}.
\end{align*}
Assume that $\gamma_1\geq\gamma_2$ and $\gamma_2=(\gamma_1+1/2)/\rho-1/2$. Then, Bai and Taqqu \cite{bt17} showed that $Z_{\gamma_1,\gamma_2}(1)\rightarrow_d Y_\rho$ as $\gamma_1\rightarrow-1/2$. (Note that if $\gamma_1\rightarrow-1/2$, then $\gamma_2\rightarrow-1/2$.)

By the representation (\ref{gamrep}), we have that  $Y_\rho\sim\mathrm{VG}_c(1,(a_\rho-b_\rho)/\sqrt{2},\sqrt{2a_\rho b_\rho})$. It was shown by Arras et al.\ \cite{aaps17} that, for any $i\geq2$, $\kappa_i(Z_{\gamma_1,\gamma_2}(1))=\kappa_i(Y_\rho)+O(-\gamma_1-1/2)$, as $\gamma_1\rightarrow-1/2$,
and combining with (\ref{optvg0}) and (\ref{optvg}) yields that, as $\gamma_1\rightarrow-1/2$, 
\begin{align*}d_{\mathrm{W}}(Z_{\gamma_1,\gamma_2}(1),Y_\rho)&\leq C_\rho\sqrt{-\gamma_1-\frac{1}{2}},\quad
d_{2}(Z_{\gamma_1,\gamma_2}(1),Y_\rho)\leq C_\rho'\Big(-\gamma_1-\frac{1}{2}\Big).
\end{align*}
where $C_\rho$ and $C_\rho'$ only depend on $\rho$. 
By employing a standard smoothing technique to the above $d_2$ metric bound, Gaunt and Li \cite{gl22} deduced the Kolmogorov distance bound
\[\sup_{z\in\mathbb{R}}|\mathbb{P}(Z_{\gamma_1,\gamma_2}(1)\leq z)-\mathbb{P}(Y_\rho\leq z)|\leq C_\rho''\Big(-\gamma_1-\frac{1}{2}\Big)^{1/3}\log\bigg(\frac{1}{-\gamma_1-1/2}\bigg), \quad \text{as $\gamma_1\rightarrow-\tfrac{1}{2}$} .\]

\appendix
 
\section{The modified Bessel function of the second kind}\label{appa}

In this appendix, we state some basic properties of the modified Bessel function of the second kind
 that are used in this review. 
Unless otherwise stated, these properties are given in Olver et al.\ \cite{olver}.

The modified Bessel function of the second kind $K_\nu(x)$ is defined for $\nu\in\mathbb{R}$ by
\[K_\nu(x)=\int_0^\infty \mathrm{e}^{-x\cosh(t)}\cosh(\nu t)\,\mathrm{d}t,\quad x > 0.
\]
For $x>0$, the function $K_\nu(x)$ is positive for all $\nu\in\mathbb{R}$.
For $\nu=m+1/2$, $m=0,1,2,\ldots$, 
\begin{equation}\label{special} K_{m+1/2}(x)=\sqrt{\frac{\pi}{2x}}\sum_{j=0}^m\frac{(m+j)!}{(m-j)!j!}(2x)^{-j}\mathrm{e}^{-x}.
\end{equation}
The modified Bessel function $K_\nu(x)$ has the following asymptotic behaviour:
\begin{eqnarray}\label{Ktend0}K_{\nu} (x) &\sim& \begin{cases} 2^{|\nu| -1} \Gamma (|\nu|) x^{-|\nu|}, & \quad x \downarrow 0, \: \nu \not= 0, \\
-\log x, & \quad x \downarrow 0, \: \nu = 0, \end{cases} \\
\label{Ktendinfinity} K_{\nu} (x) &\sim& \sqrt{\frac{\pi}{2x}} \mathrm{e}^{-x}, \quad x \rightarrow \infty,\: \nu\in\mathbb{R}.
\end{eqnarray}
The following differentiation formula holds:
\begin{align}
\label{ddbk}\frac{\mathrm{d}}{\mathrm{d}x}\big(x^\nu K_\nu(x)\big)&=-x^{\nu} K_{\nu-1}(x).
\end{align}
An integral formula of Gradshetyn and Ryzhik \cite{gradshetyn} states that for $a,\alpha>0$ and $\nu>-1/2$,
\begin{equation}\label{besint}\int_0^a t^{\nu}K_{\nu}(\alpha t)\,\mathrm{d}t =\frac{\sqrt{\pi}2^{\nu-1}\Gamma(\nu+1/2)}{\alpha^\nu}a\big(K_{\nu}(\alpha a)\mathbf{L}_{\nu-1}(\alpha a)+K_{\nu-1}(\alpha a)\mathbf{L}_{\nu}(\alpha a)\big),
\end{equation}
It is immediate from the PDF (\ref{vgpam3}) that, for $\nu>-1/2$ and $0\leq|\beta|<\alpha$,
\begin{equation}\label{pdfintegral}\int_{-\infty}^\infty |t|^{\nu}\mathrm{e}^{\beta t}K_\nu(\alpha |t|)\,\mathrm{d}t=\frac{\sqrt{\pi}(2\alpha)^\nu\Gamma(\nu+1/2)}{(\alpha^2-\beta^2)^{\nu+1/2}}.
\end{equation}
The following asymptotic approximation is readily deduced from a rescaling of the limiting form (2.13) in Gaunt \cite{gaunt ineq3}. For $\beta<\alpha$ and $\nu>-1/2$,
\begin{equation}\label{kintap}\int_x^\infty\mathrm{e}^{\beta t}t^\nu K_{\nu}(\alpha t)\,\mathrm{d}t\sim\sqrt{\frac{\pi}{2\alpha}}\frac{1}{(\alpha-\beta)}x^{\nu-1/2}\mathrm{e}^{-(\alpha-\beta)x}, \quad x\rightarrow\infty.
\end{equation}
The ratio $K_{\nu-1}(x)/K_\nu(x)$ satisfies the following bounds.  For $x>0$,
\begin{align}\label{seg1}\frac{x}{\nu-1/2+\sqrt{(\nu-1/2)^2+x^2}}<\frac{K_{\nu-1}(x)}{K_\nu(x)}&<\frac{x}{\nu-1+\sqrt{(\nu-1)^2+x^2}}, \quad \nu>1/2; \\
\label{seg3}\frac{K_{\nu-1}(x)}{K_\nu(x)}&\leq\frac{x}{\nu-1/2+\sqrt{(\nu-3/2)^2+x^2}}, \quad \nu\geq3/2,
\end{align}
with equality in (\ref{seg3}) if and only if $\nu=3/2$. The lower and upper bounds in (\ref{seg1}) were derived by Segura \cite{segura} and 
Ruiz-Antol\'{i}n and Segura \cite{rs16}, respectively. Inequality (\ref{seg3}) is given in Gaunt and Merkle \cite{gm21}.




\end{document}